\documentclass[12pt]{amsart}

\usepackage{tasks}
\usepackage{hyperref}
\usepackage{blindtext}
\usepackage{setspace}
\usepackage[utf8]{inputenc}

\usepackage{float}
\usepackage{graphicx}
 
\usepackage{subcaption}
\usepackage[utf8]{inputenc}
\usepackage[export]{adjustbox}
\usepackage{wrapfig}
\usepackage{marginnote}
\usepackage{datetime}
\usepackage{color}
\usepackage{amsrefs}
\usepackage{amsmath}
\usepackage{mathtools}
\usepackage{graphicx}  
\usepackage[top=3.5cm, bottom=2.5cm, outer=2.5cm, inner=2.5cm]{geometry}% margines
\usepackage{hyperref}
\hypersetup{linkcolor=red,
citecolor=black
}
\usepackage{amsmath}
\usepackage{float}    
\usepackage{placeins}
\newtheorem{thm}{Theorem}[section]

\newtheorem{prop}[thm]{Proposition}

\theoremstyle{definition}

\newtheorem{Alg}[thm]{Algorithm}

\usepackage{hyperref}
\usepackage{hyperref}
\hypersetup{
	colorlinks = true,
	linkcolor = {red},
	citecolor = {blue}
}

\theoremstyle{remark}
\newtheorem{re}[thm]{Remark}

\newcommand{\RNum}[1]{\uppercase\expandafter{\romannumeral #1\relax}}

\numberwithin{equation}{section}

\usepackage[utf8]{inputenc}
\usepackage[english]{babel}
\usepackage{paralist}
\usepackage{amsthm}
\usepackage{marginnote}
\usepackage{todonotes}
\usepackage{dirtytalk}
\usepackage{amssymb}

\theoremstyle{definition}

\theoremstyle{remark}

 \usepackage{bm,amsmath}

\newcommand{\tran}{^{\mathstrut\scriptscriptstyle\top}} 
 
\newcommand{\RN}[1]{%
  \textup{\uppercase\expandafter{\romannumeral#1}}%
}

\newcommand{\Rey}{\mathcal{R}e }

\newcommand{\bfzero}{\mathbf{0}}
\newcommand{\bff}{\mathbf{f}}
\newcommand{\bfu}{\mathbf{u}}
\newcommand{\bfe}{\mathbf{e}}
\newcommand{\bfv}{\mathbf{v}}

\newcommand{\bfX}{\mathbf{X}}

\newcommand{\bfV}{\mathbf{V}}
\newcommand{\bfH}{\mathbf{H}}
\newcommand{\bfx}{\mathbf{x}}

\newenvironment{Alirev}{\color{black}}{\color{black}}
\newcommand{\AAA}{\begin{Alirev}}
\newcommand{\PPP}{\end{Alirev}}

\allowdisplaybreaks[4]
\usepackage[utf8]{inputenc}

\title[Assimilating Local Data for Global Accuracy]{Global Recovery from Local Data: Interior Nudging for 2D Navier-Stokes equations in a Physical Domain}

\author[Fang]{Rui Fang}

\author[Pakzad]{Ali Pakzad$^*$}
\thanks{$^*$Corresponding author. Email: \texttt{pakzad@csun.edu}}

\address[R. Fang]{Department of Mathematics, The Ohio State University}
\email{\url{fang.1211@osu.edu}}

\address[A. Pakzad]{Department of Mathematics, California State University Northridge}
\email{\url{pakzad@csun.edu}}


\begin{document}

\begin{abstract}
In many real-world applications of data assimilation, the strategic placement of observers \AAA and sensors \PPP   is crucial for effective and efficient forecasting. Motivated by practical constraints in sensor deployment, we show that global recovery of the flow field can be achieved using observations available only in a subregion of the domain, possibly far from the boundary. We focus on the two-dimensional incompressible Navier--Stokes equations posed in a bounded physical domain with Dirichlet boundary conditions.  Building on the continuous data assimilation framework of Azouani, Olson, and Titi (2014), we rigorously prove that the assimilated solution converges globally to the true solution under suitable conditions on the nudging parameter, spatial resolution, and the geometry of the observation region, specifically, when the maximum distance from any point in the domain to the observational subregion is bounded by a constant multiple of \( \nu^{1/2} \) (in terms of scaling).  Our computational results, obtained using finite element methods over complex geometries \AAA  with both noise-free and noisy observations, \PPP support the theoretical findings and reveal even greater robustness in practice. Specifically, synchronization with the true solution is achieved even when the observational subregion lies farther from the rest of the domain than the theoretical threshold permits. Across all three tested scenarios, the local nudging algorithm performs comparably to full-domain assimilation, reaching global accuracy up to machine precision. Interestingly, observational data near the boundary are found to be largely uninformative.  This demonstrates that full observability is not necessary: carefully chosen interior observations, even far from the boundary, can suffice. 
\end{abstract}

\maketitle

\section{Introduction}\label{Section; Intro}
Data Assimilation (DA) refers to a class of methodologies that combine observational data with mathematical models to improve the accuracy of forecasts \cites{Asch-Data2016,Law-AMathematical2015}. As the volume and quality of observational data grow rapidly, thanks to advances in sensing technologies and computational power, the role and impact of DA techniques have become increasingly vital for scientific prediction. In weather and turbulence  forecasting,  observations are collected by instruments such as weather stations, sensors, or ocean buoys, all of which require significant manpower, financial resources, and logistical planning. Optimal placement of these observers is a central challenge in the design, prediction, estimation, and control of high-dimensional systems (such as weather prediction models) \cite{Chen2025}. For example, it is much easier to measure fluid velocity or temperature at shallow ocean depths than at extreme depths. Likewise, deploying a dense array of uniformly placed sensors is plainly infeasible in settings like modeling the solar wind. Therefore, from a practical standpoint, sensor placement strategies must account for resource limitations while still achieving high prediction accuracy. By identifying where data are most informative, we can reduce the number of sensors needed while improving model performance. There has been tremendous effort in recent years within the fields of control theory and inverse problems to provide guidance on optimal sensor placement for improved flow prediction \cites{zaki2021flow-reconstruction, buchta2021observation, wang2021state-estimation, encinar2019logarithmic, bewley2004footprints, suzuki2017estimation}.  In fact, it is well understood that not all observation locations contribute equally to state estimation \cites{Chen2025, zaki2021flow-reconstruction}.  Despite this, with the notable exception of the work by Biswas, Bradshaw, and Jolly \cite{Biswas2021}, the majority of existing rigorous and computational results in the direction of Continuous Data Assimilation \cite{AOT14} (CDA; see Subsection~\ref{Section;CDA}) assume \emph{global observability}; that is, observational data are assumed to be uniformly available throughout the entire spatial domain. This assumption, while mathematically convenient, poses serious limitations for real-world applications. Our work challenges this framework by showing that it is possible to recover the full system dynamics using only local observational data from a subregion of the domain.

In this paper, we conduct a rigorous analysis, accompanied by a comprehensive computational study, of a \textbf{global DA algorithm that relies solely on local observations}. Our focus is on the two-dimensional Navier-Stokes equations (NSE) posed in a bounded physical domain \( \Omega \), subject to Dirichlet boundary conditions
\begin{equation} \label{NSE}
\partial_t \bfu + (\bfu \cdot \nabla) \bfu - \nu \Delta \bfu + \nabla p = \bff, \quad \nabla \cdot \bfu = 0, \quad \text{in } \Omega \times (0, \infty),
\end{equation}
with the initial condition \( \bfu_0 = \bfu(x,0) \), which is typically not fully known in practice.
 Here, \( \bfu = \bfu(x,t) \) denotes the \say{true} velocity field, \( p = p(x,t) \) is the pressure, \( \bff = \bff(x,t) \) is a given external body force, and \( \nu > 0 \) is the kinematic viscosity. When using the NSE, or any other relevant nonlinear model, for weather prediction, the exact initial condition \( \bfu_0 \) is typically not fully known a priori \cite{K03}. This uncertainty in the initial data can lead to significant challenges, as even small errors may grow exponentially in time due to the nonlinear nature of the system \cites{lorenz1963deterministic, lorenz1963predictability, lorenz1965study, K03}. To overcome this difficulty, in DA, the observational measurements are directly inserted into a mathematical model \cites{hoke1976initialization,kistler1974study}.  

Motivated by \cite{AOT14}, we now introduce the framework of \emph{local DA}, where the solution \( \bfu \) is only partially and locally known, specifically in a subregion located \emph{far} from the boundary. Specifically, we assume that coarse spatial observations of the velocity field are available \emph{only} within an interior subdomain \( \Omega_0 \subset \Omega \) (see Figure~\ref{fig:Bezier domain}), and not throughout the entire physical domain $\Omega$.  These observations are represented by an interpolation operator \( I_H(\bfu) \), which acts on \( \bfu \) to extract its values on a coarse spatial grid of resolution \( H \) supported in \( \Omega_0 \).   The corresponding local DA algorithm is formulated as the following modified Navier–Stokes system with the same Dirichlet boundary conditions in the physical domain \( \Omega \)
\begin{equation} \label{LocDAnse}
\partial_t \bfv + (\bfv \cdot \nabla) \bfv - \nu \Delta \bfv + \nabla q = \bff + \mu \big( I_H(\bfu) - I_H(\bfv) \big), \quad \nabla \cdot \bfv = 0,  \quad \text{in } \Omega \times (0, \infty), 
\end{equation}
where
\begin{center}
\textit{$I_H$: localization operator that incorporates observational data} \\
\textit{only from a designated subregion $\Omega_0$ of the domain,}
\end{center}
\vspace{0.2cm}
and the initial value of $\bfv(x,0)= \bfv_0$ is arbitrary. Here, \( \bfv = \bfv(x,t) \) is the \say{approximate} velocity field, and \( \mu > 0 \) is a relaxation (nudging) parameter. 

\subsection*{Results in this paper} Having established the framework in \eqref{LocDAnse},  we provide a rigorous analysis showing that it is possible to recover the full flow field using velocity observations collected only from an interior portion of the domain, without requiring measurements near the boundaries.  Specifically, we prove in Theorem~\ref{Theorem;Main} that if no part of the domain lies too far from the observational subregion (more precisely, if the maximum distance from any point in the domain to the observation region remains below a fixed threshold proportional to the square root of viscosity in terms of scaling\footnote{This is the same scaling that governs the thickness of the classical Prandtl boundary layer\cite{Prandtl1904_BL}, suggesting a potential  connection between the effectiveness of interior observations and the physical structure of near-wall viscous flow.}, \( \delta \lesssim \nu^{1/2} \)) then, under standard conditions on the nudging parameter and spatial resolution, the locally assimilated solution \eqref{LocDAnse} converges globally to the true solution of the $2D$ NSE\eqref{NSE}. This result demonstrates that achieving global accuracy does not require global data, and  quite surprisingly, not even data near the boundary. However, carefully chosen, localized observations from the interior of the domain can be sufficient.

We demonstrate the efficacy and robustness of the algorithm through extensive computational studies \AAA  using both noise-free and noisy observations. \PPP    In all experiments, we implement the finite element method for simulations over complex geometries with irregular domains and realistic boundary conditions. A key technical challenge, and novelty, of this work lies in constructing a local, coarse finite element mesh for the interpolation operator that is compatible with the global fine mesh used to solve the full system. This compatibility ensures mathematical consistency and numerical stability while enabling localized DA. Three flow scenarios are tested
\begin{enumerate}
    \item Flow over a flat obstacle (see Section \ref{Section; SIM1}),  
    \item  Couette/Shear flow in an annulus (see Section \ref{Section; SIM2}), 
    \item Flow driven by a body force in a disk with an off-center obstacle (see Section \ref{Section; SIM3}).  
\end{enumerate}
All flows are simulated within a bounded physical domain under Dirichlet boundary conditions. All three simulations are in agreement with the theoretical results proved in Theorem~\ref{Theorem;Main}. We observe that the local nudging algorithm achieves global  synchronization within machine precision, typically at a slightly (but not significantly) slower rate compared to full-domain nudging,  when the data are collected from an interior subregion far from the boundaries. In nearly all simulations, observational data near the boundaries are found to be uninformative for the nudging algorithm. To ensure that this effect was not simply due to variations in data density or the number of observational points, we also varied the area of the local DA region. Nevertheless, the same trend was consistently observed: interior observations led to effective synchronization, while boundary observations remained largely uninformative.

The theoretical strict conditions on the parameters, namely the nudging parameter \( \mu \), the density of observational points \( H \), and the distance \( \delta \) between the observational subregion and the rest of the domain, can be significantly relaxed in practice.  Even when the local  observation region lies farther from the boundaries than the theoretical threshold permits, the local assimilation algorithm still performs well, achieving global synchronization with the true solution up to machine precision. This theoretical restrictions is well accepted in the frame work of the data assimilation  (see, for instance, \cites{Rebholz_Zerfas2021,CIBIK2025117526,Gesho-Acomputational2016, JP23, CAO2022103659, Cao-Algebraic2021, GLRVZ21,PhysRevFluids.9.054602}). 

\subsection{\textit{Global} Continuous Data Assimilation}\label{Section;CDA}

Using observational data to improve the predictive capability of mathematical models dates back to the work of Rudolph Kalman in 1960 \cite{Kalman1960}, and of Charney, Halem, and Jastrow in 1969 \cite{Charney-Use1969}, following the advent of satellite-borne observation systems that began producing climate data. In the same direction, DA aims to recover the fine-scale dynamics of a system using only coarse, partial observations \cites{K03, Daley1991}. A notable contribution in this area is the work of Azouani, Olson, and Titi \cite{AOT14}, who introduced a nudging-based DA algorithm (a.k.a. CDA) that is mathematically rooted in earlier work of Foias and Prodi \cite{Foias-Sur1967} on determining functionals. One of the key advantages of their approach lies in its simplicity of implementation and its compatibility with rigorous mathematical analysis, making it a compelling alternative to other data-driven techniques.

To describe the CDA algorithm, suppose $\bfu(t)$ is a solution of a dissipative dynamical system governed by
\begin{equation}
    \frac{d\bfu}{dt} = F(\bfu),
\end{equation}
where the initial data $\bfu(0)=\bfu_0$ is not fully known. We denote the observable state of the solution by \( I_h(\bfu) \), where \( I_h \) is an autonomous, bounded linear operator representing an orthogonal projection onto a finite-dimensional subspace of \( L^2(\Omega) \). This operator satisfies the approximation property
\begin{equation}\label{LocI_h}
\|\phi - I_h(\phi)\|_{L^2(\Omega)} \leq C_1 h \|\nabla \phi\|_{L^2(\Omega)},
\end{equation}
for all \( \phi \in H^1(\Omega) \), where \( C_1 > 0 \) is a constant independent of the discretization parameter \( h \). The operator \( I_h \) approximates the true state by utilizing only the information available from  global coarse spatial observations $h$. In the context of incompressible flows considered in this work, \( I_h \) is understood to first interpolate the observed data, and then apply an orthogonal projection onto the space of divergence-free (solenoidal) vector fields. Consequently, \( I_h \) is idempotent and self-adjoint, i.e., \( I_h^2 = I_h \) and \( I_h^* = I_h \).  Examples include the orthogonal projection onto low Fourier modes \cite{JonesTiti}  and, as a more physically motivated example, interpolation using piecewise polynomial finite element spaces \cite{JP23}.

Using the CDA algorithm introduced in \cite{AOT14}, the observations $I_h(\bfu(t))$ are incorporated into an auxiliary system given by
\begin{equation}
    \frac{d\bfv}{dt} = F(\bfv) + \mu \left( I_h(\bfu) - I_h(\bfv) \right),
\end{equation}
with an \textit{arbitrary} initial condition \( \bfv_0 \). Under appropriate conditions, that is, sufficiently large \( \mu \) and sufficiently small \( h \), one can have
$$ \bfv(t) \to  \bfu(t) $$  
exponentially fast in a suitable function space. 

Interpolation-based nudging has been extensively studied from both theoretical and computational perspectives. Foundational theoretical results include the convergence and error analysis for the $2D$ NSE by Azouani, Olson, and Titi \cite{AOT14}, later extended to the $3D$ setting by Biswas and Price \cite{Biswas_3D}. Additional error estimates using the $L^2$ projection operator were developed by Rebholz and Zerfas \cite{Rebholz_Zerfas2021}. Interpolation-based nudging has a wide application to diverse fluid dynamics problems, as seen in \cites{CAO2022103659, FJT, Farhat-Data2020, Jolly-Continuous2019, Pei-Continuous2019, Bessaih-Continuous2015, Leo2022, FLT316, Jolly-Adata2017, ALT16, FLT, Gesho-Acomputational2016, GLRVZ21, newey2024model, cibik2025data}. Beyond idealized models, recent research has demonstrated that data assimilation via interpolant nudging remains effective even with model errors or noisy observations \cites{hawkins2024removing, cibik2025data, Larios_Pei2020, chen2023data}.   Rebholz and Zerfas~\cite{LeoAlgebricDA} introduced an efficient algebraic nudging operator that enables continuous data assimilation to be incorporated into finite element and legacy codes through a simple diagonal matrix,  while preserving long-time stability and synchronization.  The selection of the nudging parameter $\mu$, which controls assimilation strength, has been shown to be critical for performance, and several works have proposed adaptive or optimized selection methods \cites{CIBIK2025117526, Diegel2025}. Alternative data assimilation frameworks such as ensemble Kalman filters and $3D/4D$-VAR methods \cites{K03, Blomker2013}, but interpolation-based nudging remains attractive for its simplicity and provable convergence properties.

\subsection*{Organization of this paper} In Section \ref{Section; Prelim}, we introduce the inequalities and preliminary results on the NSE used in our analysis. Section \ref{Section; FEM} provides a  brief background on the finite element method and the numerical schemes used throughout our simulations; more importantly, we detail how the global Direct Numerical Simulation (DNS) mesh is constructed from the coarse local observational mesh. In Section \ref{Section; Analysis}, we state and prove our main theoretical results. The three numerical experiments, each demonstrating and extending beyond the analytical results, are described in Sections \ref{Section; SIM1}, \ref{Section; SIM2}, and \ref{Section; SIM3}. Finally, in Section \ref{Section; Discussion}, we discuss potential improvements and future directions.

\section{Preliminaries}\label{Section; Prelim}

The  results in this section are standard and can be found in any classical text on the NSE \cite{FMRT01} or functional analysis \cite{Brezis2011}.

Let $\bfx=(x_1, x_2) \in  \mathbb{R}^2$ and  \( \Omega \subset \mathbb{R}^2 \) be a bounded domain with \( C^2 \) boundary. Let \( 1 \leq p \leq \infty \), and denote by \( p' \) the conjugate exponent defined by \( \frac{1}{p} + \frac{1}{p'} = 1 \). Assume that \( \phi \in L^p(\Omega) \) and \( \psi \in L^{p'}(\Omega) \). Then H\"older's inequality states that
\begin{equation} \label{Holder}
\|\phi \, \psi\|_{L^1(\Omega)} \leq \|\phi\|_{L^p(\Omega)} \, \|\psi\|_{L^{p'}(\Omega)}.
\end{equation}
Additionally, for any \( a, b \geq 0 \) and \( \lambda > 0 \), Young's inequality takes the form
\begin{equation} \label{Young}
a b \leq \lambda a^p + \frac{1}{p'} \left( p \lambda \right)^{-\frac{p'}{p}} b^{p'}.
\end{equation}
For any \( \phi \in H^1(\Omega) \), Ladyzhenskaya's inequality states
\begin{equation}\label{Lady-ineq}
    \|\phi\|_{L^4(\Omega)} \leq 2^{1/4} \|\phi\|^{1/2}_{L^2(\Omega)} \|\nabla \phi\|^{1/2}_{L^2(\Omega)}.
\end{equation}
Let \( d(\mathbf{x}) = \operatorname{dist}(\mathbf{x}, \partial \Omega) \). Then, for any \( \mathbf{u} \in W_0^{1,p}(\Omega) \), Hardy's inequality \cite{Brezis2011} in p.~313 states that there exists a dimensionless constant \( C_p \) such that
\begin{equation}\label{Hardy}
    \left\| \frac{\mathbf{u}}{d} \right\|_{L^p(\Omega)} \leq C_p \| \nabla \mathbf{u} \|_{L^p(\Omega)}.
\end{equation}

To support the forthcoming analysis, we first recall some fundamental results  of \eqref{NSE}. These results are well-known and documented in the literature (see, for example, \cites{T77, FMRT01}). Denote by \( \lambda_1 \) the smallest eigenvalue of the Stokes operator, then for $\phi \in H^1$,  the Poincar\'e's  inequality reads as 
\begin{equation}\label{Poincare}
\lambda_1 \|\phi\|^2 \leq  \|\nabla \phi\|^2.
\end{equation}
A key non-dimensional parameter that captures the complexity of the system is the Grashof number, given by
\begin{equation}\label{Grashof}
G \coloneqq \frac{1}{\nu^2 \lambda_1} \|\bff\|_{L^2(\Omega)}.
\end{equation}
\begin{thm}[Existence and Uniqueness of Strong Solutions in $2D$]
Suppose \( \bfu_0 \in H^1(\Omega) \) and \( \bff \in L^{\infty}\big((0, \infty); L^2(\Omega)\big) \). Then, for any \( T > 0 \), the initial value problem \eqref{NSE} admits a unique strong solution \( \bfu \) satisfying
\[
\bfu \in C\big([0,T]; L^2(\Omega)\big) \cap L^2\big(0,T; H^1(\Omega)\big), \quad \partial_t \bfu \in L^2\big(0,T; L^2(\Omega)\big).
\]
\end{thm}

\begin{thm}\label{NSEBounds}
Fix \( T > 0 \), and let \( G \) be the Grashof number as defined in \eqref{Grashof}. Suppose that \( \bfu \) is the solution of \eqref{NSE} corresponding to the initial value \( \bfu_0 \). Then, there exists a time \( t_0 \), depending on \( \bfu_0 \), such that for all \( t \geq t_0 \), we have:
\[
\|\bfu(t)\|^2 \leq 2 \nu^2 G^2,
\]
and
\[
\int_t^{t+T} \|\nabla \bfu(\tau)\|^2 \, d\tau \leq 2 \left(1 + T \nu \lambda_1\right) \nu G^2.
\]
\end{thm}

\begin{prop}[Uniform Gr\"onwall Inequality]\label{Gronwall} 
Let $T>0$ be fixed. Suppose
$$\frac{d}{dt} Y + \alpha(t) Y \leq 0,$$
where 
$$\limsup_{t \to \infty} \int_t^{t+T} \alpha(s) ds \geq \alpha_0 >  0,$$
then $Y(t) \to 0$, exponentially, as $t \to \infty$.
\end{prop}

\section{Mesh Construction and Numerical Implementation}\label{Section; FEM}
In this section, we introduce the necessary notation for the finite element method and provide details on how the global mesh used for solving the full system is constructed from the coarse local observational mesh, an essential component of our approach.

\begin{figure}
    \centering
    \begin{subfigure}{0.49\textwidth}
        \centering
        \includegraphics[width=\linewidth]{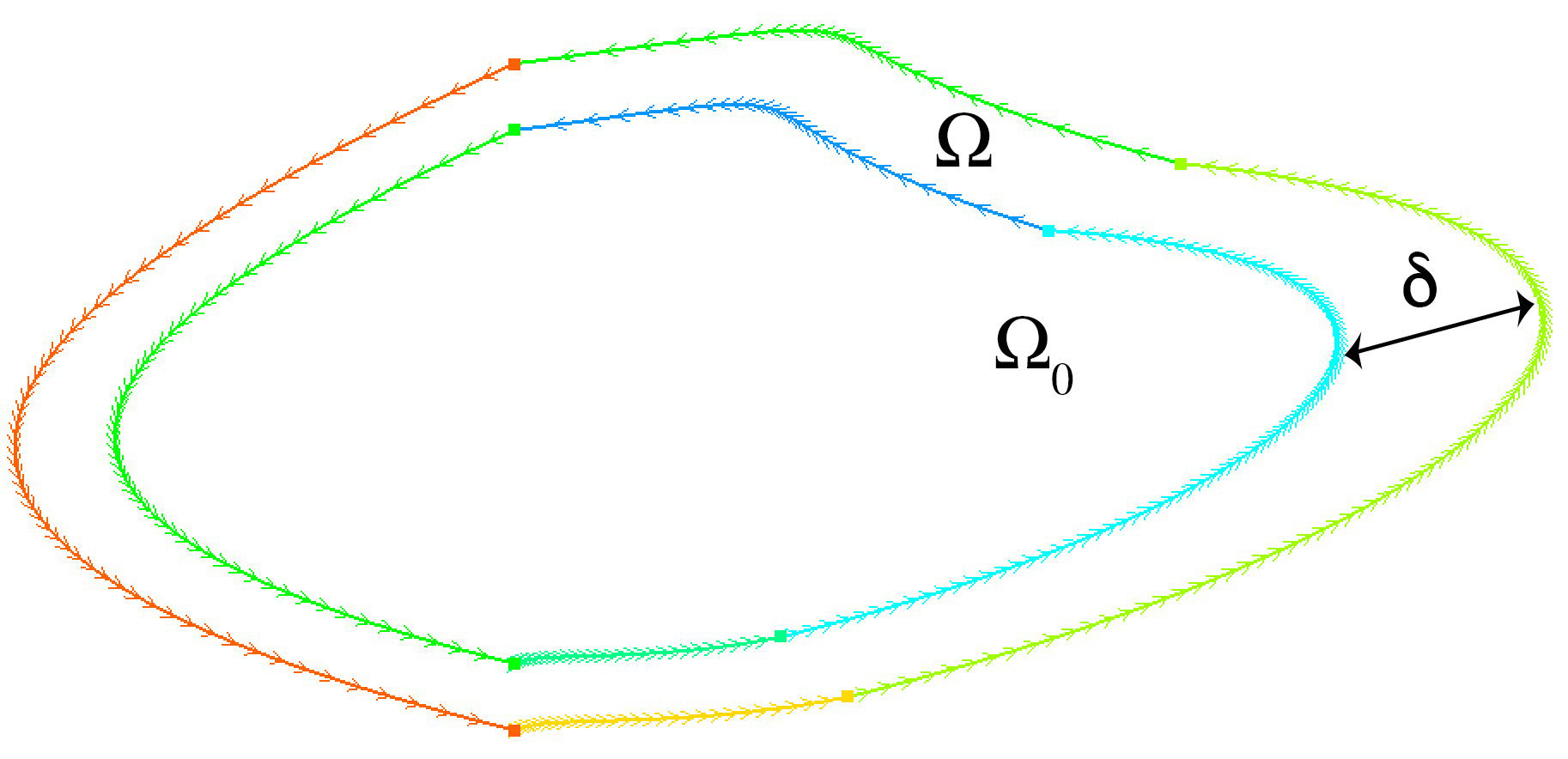}
        \caption{$\Omega$: Physical Domain. $\Omega_0$: Observation  subdomain.}
        \label{fig:Bezier domain}
    \end{subfigure}
    \begin{subfigure}{0.49\textwidth}
        \centering
        \includegraphics[width=\linewidth]{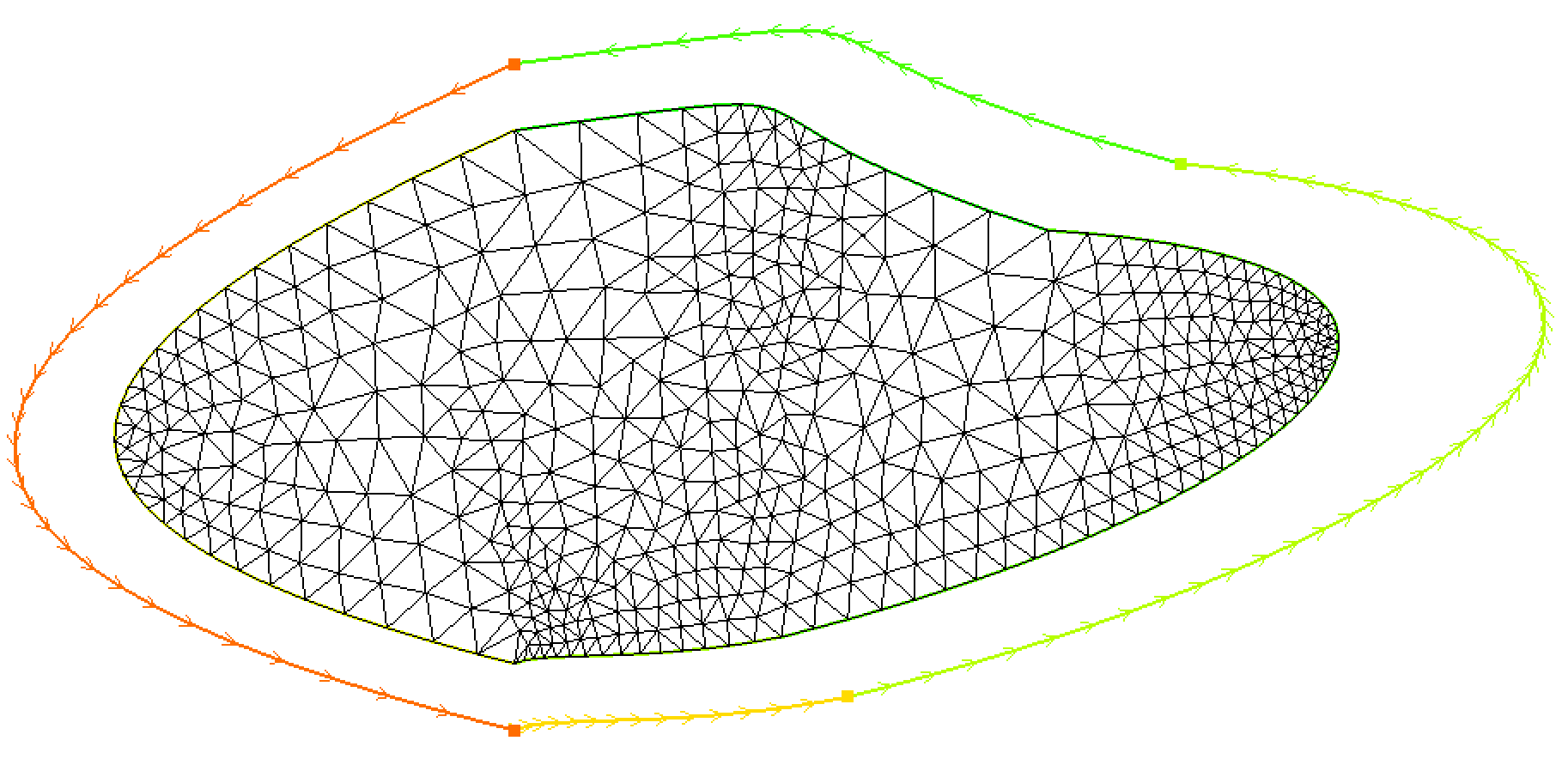}
       \caption{Local coarse spatial resolution $H$.}
        \label{fig:Bezier_domain_inner_coarse}
    \end{subfigure}
    \begin{subfigure}{0.49\textwidth}
        \centering
        \includegraphics[width=\linewidth]{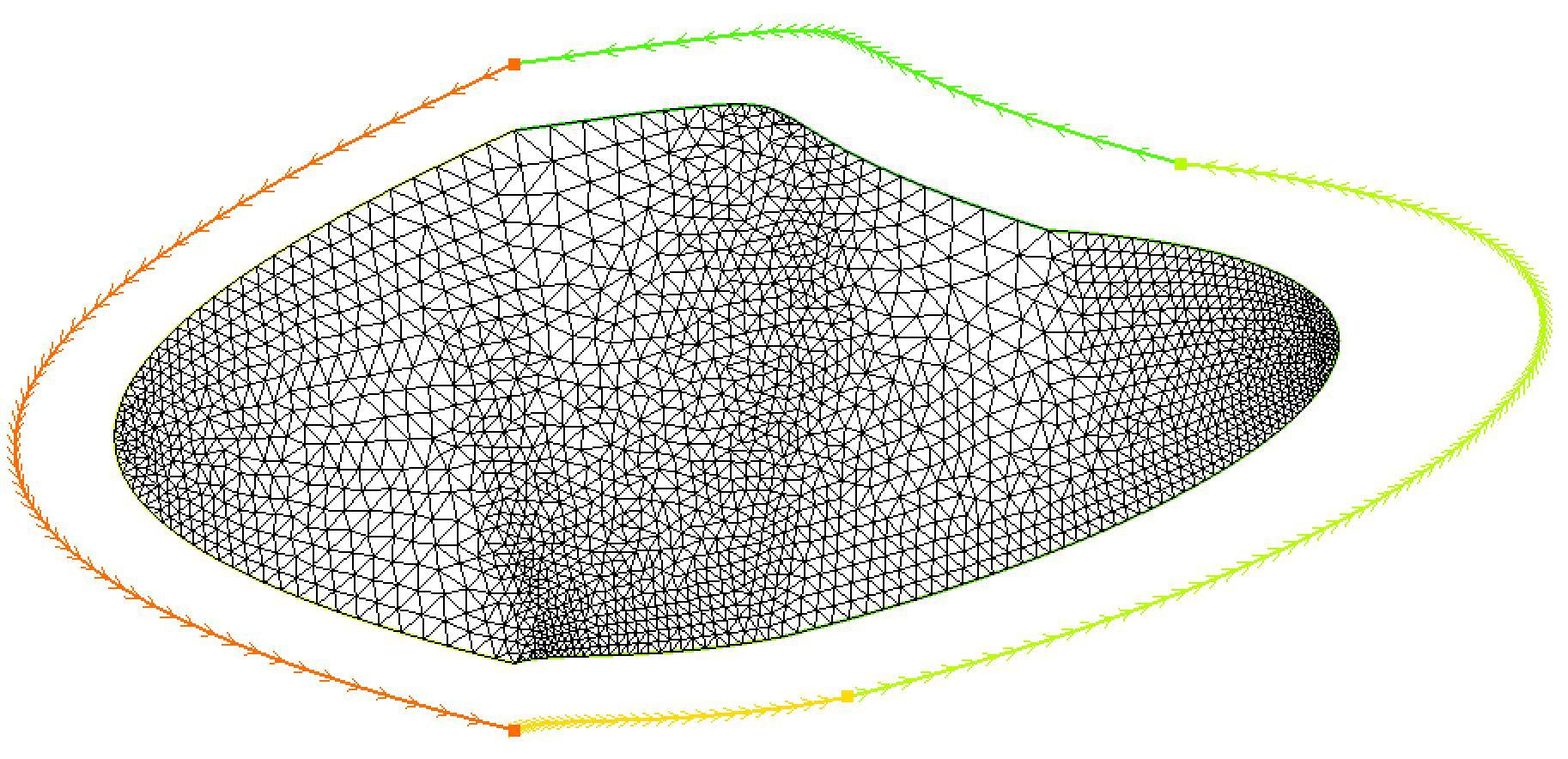}
           \caption{Local DNS mesh.}
           \label{fig:Bezier_domain_inner_fine}
    \end{subfigure}
    \begin{subfigure}{0.49\textwidth}
        \centering
        \includegraphics[width=\linewidth]{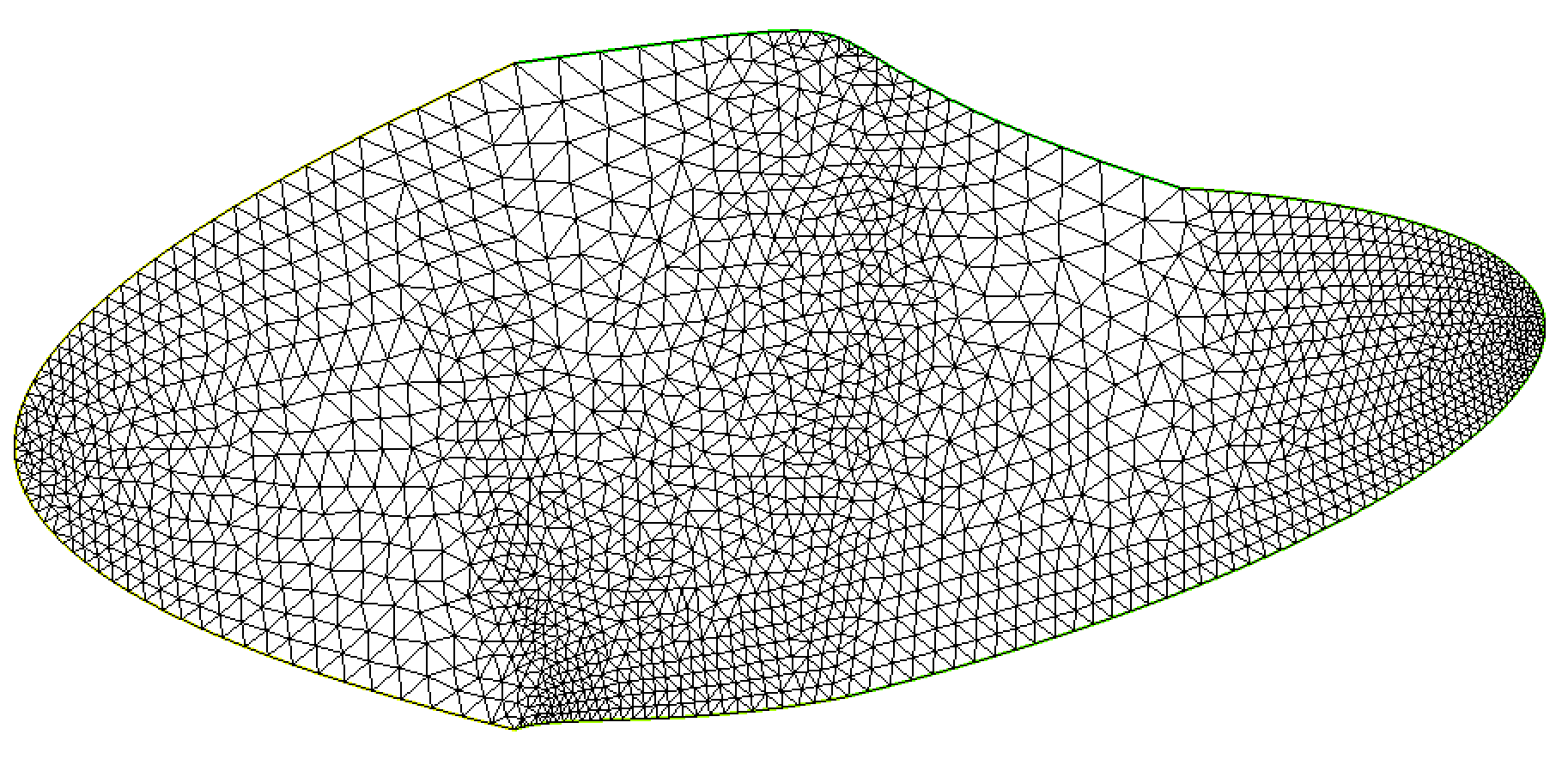}
        \caption{Global DNS mesh $h$.}
        \label{fig:Bezier_domain_fine}
    \end{subfigure}
    \caption{Domain mesh setup: from local coarse observations to global DNS mesh.}
    \label{fig:case1}
\end{figure}

Let \( \Omega \subset \mathbb{R}^2 \) be a bounded domain with \( C^2 \) boundary, and let \( \Omega_0 \subset \Omega \) be a subdomain; see Figure \ref{fig:Bezier domain}. We define \( \delta \) to be the maximum distance between any point in \( \Omega \) and the subdomain \( \Omega_0 \); that is
\[
\delta := \sup_{x \in \Omega} \operatorname{dist}(x, \Omega_0).
\]
From this point onward, we denote by \( \Omega_{\delta} \) the region between \( \Omega \) and \( \Omega_0 \); that is
\[
\Omega_{\delta} := \Omega \setminus \overline{\Omega_0}.
\]
Denote the natural function spaces for velocity and pressure, respectively, by
\begin{align*}
    \bfX & := \bfH_0^1(\Omega) = \left\{ \bfu \in \bfH^1(\Omega): \bfu = 0 \ \text{on} \ \partial \Omega \right\},\\
    Q &:= L_0^2(\Omega) = \left\{ q \in L^2(\Omega): \int_\Omega q \, dx = 0 \right\}, \\
    \bfV& := \left\{ \bfu \in \bfX: \ (\nabla \cdot \bfu, q) = 0 \ \text{for all } q \in Q \right\}.
    \end{align*}
Denote by $\Omega_h$ a regular, conforming triangulation of the domain $\Omega$, and  let $\bfX^h \subset \bfX$ and $Q^h \subset Q$ be an
inf-sup stable pair of discrete velocity--pressure spaces (see \cite{J16} for more details). Additionally, we define the discretely divergence-free space $\bfV^h$ as
\[
\bfV^h := \left\{ \bfv^h \in \bfX^h : (q^h, \nabla \cdot \bfv^h) = 0 \quad\text{for all } q^h \in Q^h \right\}.
\]
For all simulations, following  Algorithm \ref{Algorithm},  we employ the IMEX (implicit-explicit) scheme for temporal discretization to avoid solving a nonlinear problem at each time step. The spatial discretization is based on the finite element method. We assume the velocity--pressure finite element spaces \( (\bfX^h, Q^h) = (\mathbb{P}_2, \mathbb{P}_1) \), known as the Taylor--Hood elements, which satisfy the inf-sup stability condition, see \cite{L08} Chapter~9 for detailed description on the scheme. 

\begin{Alg}\label{Algorithm}
Given body force $\bff \in L^{\infty} ((0, \infty), L^2(\Omega))$,  initial condition $\bfv_0 \in \bfV^h$ , reference solution $\bfu_{n+1}$, and  $(\bfv^h_n, q_n^h) \in \bfX^h \times Q^h$,  compute $(\bfv_{n+1}^h, q_{n+1}^h) \in \bfX^h \times Q^h$ satisfying
\begin{multline}\label{IMEX-FEM}
(\frac{\bfv_{n+1}^{h} - \bfv_{n}^{h}}{\Delta t},\Theta ^{h}) + b(\bfv^h_n,\bfv_{n+1}^{h},\Theta^{h}) + \nu \, (\nabla \bfv_{n+1}^{h},\nabla \Theta^{h})  - (q_{n+1}^{h}, \nabla \cdot \Theta^{h})
\\ =  (\bff_{n+1},\Theta^{h})  - \mu \, (I_H(\bfv^h_{n+1} - \bfu_{n+1}), \Theta^h),
\end{multline}
\begin{equation}\label{incomp}
(\nabla \cdot \bfv_{n+1}^{h},r^{h}) = 0,
\end{equation}
for all $(\Theta^{h},  r^{h}) \in  \bfX^h \times Q^h.$
\end{Alg}
In order to minimize uncertainties, maintain full control over the resolution of observations, and assess the accuracy of the reconstruction, we extract the observations from an independent DNS. While \( I_H(\bfu_{n+1}) \) denotes the local observations restricted to the subdomain \( \Omega_0 \) at a coarse spatial resolution and at time step \( t_{n+1} \), the global approximating solution \( \bfv^h \), initialized from an arbitrary initial condition, is computed using Algorithm~\ref{Algorithm}. Note that at time \( t_{n+1} \), the relevant data about the true solution have already been acquired; therefore, \( I_H(\bfu_{n+1}) \) is interpreted as the most recent data gathered locally throughout $\Omega_0$. In all experiments, the algorithms are implemented using the public domain finite element software \texttt{FreeFEM++} \cite{FreeFEM}. Throughout all our numerical tests, the precision and robustness of the reconstruction are assessed by comparing the locally assimilated solution \( \bfv^h \), computed using the above algorithm, with the reference solution \( \bfu^h \) obtained from DNS.

\subsection*{Local-to-Global Mesh Construction}
One challenge in our simulation lies in constructing a local coarse finite element mesh for the interpolation operator that is compatible with the global DNS mesh used to solve the whole system. This compatibility ensures mathematical consistency and numerical stability, while enabling localized data assimilation. Before presenting the main results, we demonstrate this process on a general two-dimensional domain. Consider a B\'ezier surface representing the bounded physical domain \( \Omega \subset \mathbb{R}^2 \), and let \( \Omega_0 \subset \Omega \) denote the observation subdomain where local coarse spatial measurements are available;  Figure \ref{fig:Bezier domain}. We first construct the coarsest mesh of size \( H \) within \( \Omega_0 \), corresponding to the resolution at which local observational data are collected, see Figure \ref{fig:Bezier_domain_inner_coarse}. Each coarse grid cell is then split into more  sub-triangles to generate a fine mesh, consistent with the resolution of DNS, within the subdomain \( \Omega_0 \), see Figure \ref{fig:Bezier_domain_inner_fine}. Finally, this fine mesh is extended to cover the entire domain \( \Omega \), with \( h \) denoting the finest global DNS mesh size;  Figure \ref{fig:Bezier_domain_fine}.

\section{Rigorous Analysis;  Synchronization}\label{Section; Analysis}

\begin{thm}\label{Theorem;Main}
Let \( \Omega \subset \mathbb{R}^2 \) be a bounded domain with \( C^2 \) boundary, and let \( \Omega_0 \subset \Omega \) be a subdomain with 
$
\delta := \sup_{x \in \Omega} \operatorname{dist}(x, \Omega_0).
$ 
Let \( \bfu \) be a solution of the incompressible two-dimensional NSE \eqref{NSE} on \( \Omega \), equipped with no-slip Dirichlet boundary conditions. Let \( \bfv \) be the approximating solution given by \eqref{LocDAnse}, where \( I_H \) is the localization operator that incorporates observational data from the subregion \( \Omega_0 \), and satisfies \eqref{LocI_h}. Then,
\[
\|\bfu - \bfv\|_{L^2(\Omega)} \to 0,
\]
provided that the following conditions hold:
\[
\AAA 4\, C^2_1 \mu H^2 \leq \nu \PPP, \quad \mu \geq 4\, G^2 \nu \lambda_1, \quad \text{and} \quad \delta \leq \frac{C_1}{C_p} \, H.
\]
\end{thm}

\begin{re}\textbf{The boundary layer condition  can be relaxed throughout the practice/simulation}.
The above restrictions on the parameters \( \mu \), \( H \), and particularly the condition \( \delta \lesssim H \sim \mathcal{O}( \nu^{1/2}) \), may appear overly conservative from a practical standpoint. However, such limitations are well understood in the context of data assimilation and its rigorous analytical framework (see, for instance, \cites{Rebholz_Zerfas2021,CIBIK2025117526,Gesho-Acomputational2016, JP23, CAO2022103659, Cao-Algebraic2021, GLRVZ21,PhysRevFluids.9.054602}), where the bounds arise from classical inequalities and energy estimates that are not easily sharpened, but are nonetheless standard and widely accepted.  These theoretical restrictions are significantly relaxed in practice, as demonstrated throughout our simulations. Consequently, developing less restrictive analytical criteria, possibly through alternative mathematical techniques or problem-specific refinements, that better reflect the robustness observed in computational experiments remains an open and compelling direction for future research.
\end{re}

\begin{proof}
 Subtracting \eqref{NSE} and \eqref{LocDAnse}, the difference $\bfe= \bfu-\bfv$ satisfies 
 $$\partial_t \bfe + \bfe \cdot \nabla \bfu + \bfv \cdot \nabla \bfe - \nu \Delta \bfe +\nabla (p-q)= - \mu I_H (\bfe), \quad \nabla \cdot \bfe = 0.$$
Using incompressibility, and after $L^2$ multiplication of the above difference equation with $\bfe$, we arrive at
\begin{equation}\label{Energy1}
 \frac{1}{2}\frac{d}{dt} \|\bfe\|_{L^2(\Omega)}^2 + \nu \|\nabla \bfe\|_{L^2(\Omega)}^2 = - (\bfe \cdot \nabla \bfu, \bfe) -  \mu \|I_H(\bfe)\|^2_{L^2(\Omega_0)}:= \RN{1} + \RN{2}.   
\end{equation}

Using H\"older's and Young inequalities \eqref{Holder} and \eqref{Young}, along with  \eqref{Lady-ineq}, Term $\RN{1}$ in the above is estimated as 
\begin{equation}\label{Energy1-1}
\begin{split}
 |(\bfe \cdot \nabla \bfu, \bfe)| \leq  \|\bfe\|^2_{L^4(\Omega)} \| \nabla \bfu\|_{L^2(\Omega)} &\leq 2^{\frac{1}{2}}\| \bfe\|_{L^2(\Omega)} \| \nabla \bfe\|_{L^2(\Omega)} \|\nabla \bfu\|_{L^2(\Omega)}\\
 & \leq \frac{\nu}{2} \| \nabla \bfe\|^2_{L^2(\Omega)} + \frac{1}{\nu} \|\nabla \bfu\|^2_{L^2(\Omega)}  \| \bfe\|^2_{L^2(\Omega)}. 
\end{split}
\end{equation}
Next we estimate term $\RN{2}$ in \eqref{Energy1}. First note that 
\begin{equation}\label{Energy1-3}
   \| \bfe\|^2_{L^2(\Omega)}= \| \bfe\|^2_{L^2(\Omega_0)}+ \| \bfe\|^2_{L^2(\Omega_\delta)} \leq  \| I_H(\bfe)\|^2_{L^2(\Omega_0)}+ \|(I-I_H)(\bfe)\|^2_{L^2(\Omega_0)} + \| \bfe\|^2_{L^2(\Omega_\delta)}, 
\end{equation}
The second term on the right-hand side of the above inequality can be bounded as follows
 $$ \|(I-I_H)(\bfe)\|^2_{L^2(\Omega_0)} \leq (C_1H)^2 \|\nabla\bfe\|^2_{L^2(\Omega_0)} \leq (C_1H)^2 \|\nabla\bfe\|^2_{L^2(\Omega)}.$$
Then, using Hardy's inequality~\eqref{Hardy}, we obtain the following estimate for the last term in~\eqref{Energy1-3}
\[
\| \mathbf{e} \|^2_{L^2(\Omega_\delta)} \leq (C_p \delta)^2 \| \nabla \mathbf{e} \|^2_{L^2(\Omega_\delta)} \leq (C_p \delta)^2 \| \nabla \mathbf{e} \|^2_{L^2(\Omega)}.
\]

Now from \eqref{Energy1-3}, and the above two estimates, one can obtain  
\begin{equation}\label{Energy1-4}
    \begin{split}
    \| I_H(\bfe)\|^2_{L^2(\Omega_0)} \geq \|\bfe\|^2_{L^2(\Omega)}  -  2 \, \max\{C^2_1 H^2, C^2_p \delta^2 \}  \|\nabla\bfe\|^2_{L^2(\Omega)}.
    \end{split}
\end{equation}
With assumptions $\delta \leq \frac{C_1}{C_p} \, H$ and \AAA $4\, C^2_1 \mu H^2 \leq \nu $ \PPP,  and using estimates \eqref{Energy1-1} and \eqref{Energy1-4} in \eqref{Energy1}, we get the following
$$\frac{d}{dt} \|\bfe\|_{L^2(\Omega)}^2 + \AAA \left[ 2 \mu - \frac{2}{\nu}\|\nabla \bfu\|_{L^2(\Omega)}^2\right]\,\PPP  \,  \|\bfe\|_{L^2(\Omega)}^2 \leq 0.$$
Next, using the bounds established in Theorem~\ref{NSEBounds}, one can choose \( T = \frac{1}{\nu \lambda_1} \) and \( \mu \geq 4 \, G^2 \nu \lambda_1 \) to obtain
\[
\limsup_{t \to \infty} \int_t^{t+T} \left(2 \mu - \frac{2}{ \nu} \|\nabla \bfu\|_{L^2(\Omega)}^2\right) \, ds \geq \alpha_0=: 2 G^2  > 0.
\]
Then, applying the standard Gr\"onwall inequality (Proposition~\ref{Gronwall}) with 
\[
\alpha(t) = 2 \mu - \frac{2}{\nu} \|\nabla \bfu\|_{L^2(\Omega)}^2,
\]
it follows that \( \|\bfe\|_{L^2(\Omega)}^2 \to 0 \) exponentially as \( t \to \infty \).
\end{proof}

\begin{re}
We note that, under higher regularity assumptions such as $\nabla u \in L^{\infty}(\Omega)$, as considered by García-Archilla, Novo, and Titi \cite{GNT_2020}, there may be potential to improve the upper bound on the nudging parameter $\mu$. A detailed investigation of this possibility is left to the interested reader.

\end{re}

\section{Computational Study I: Flow Over a Flat Obstacle}\label{Section; SIM1}

In this section, using Algorithm \ref{Algorithm},  we conduct a numerical test of flow over a flat obstacle, which is based on the setup described in \cites{GLRVZ21, CIBIK2025117526}. The computational domain is a rectangular channel, $[-7,20] \times [0,20]$, that contains a flat plate obstacle of size $0.125 \times 1$, centered at $(0.075,0)$. \AAA  To study the effect of observation placement on local data assimilation, we
partition the computational domain $\Omega$ into three physically distinct
subdomains:
\[
\Omega
=
\Omega_{\mathrm{obst}}
\cup
\Omega_{\mathrm{bulk}}
\cup
\Omega_{\mathrm{wall}},
\]
whose interiors are pairwise disjoint. These subdomains are defined as follows:

\begin{itemize}
    \item  $\Omega_{\text{obst}}: $ Obstacle-adjacent Region, surrounding the  immersed obstacle,  with $2l$ thickness extending beyond both the length and width of the plate.
    
    \item  \(\Omega_{\text{bulk}}\): Interior bulk Region, located away from all the physical boundaries.
    
    \item \(\Omega_{\text{wall}}\):  Wall-adjacent Region, lying near the outer boundary or channel walls,  with thickness $l$ near the top and bottom walls.
\end{itemize}
The union of the two boundary-adjacent regions is denoted by
\[
\Omega_{\mathrm{bdry}}
:=
\Omega_{\mathrm{obst}} \cup \Omega_{\mathrm{wall}}.
\]

\PPP 

\begin{figure}[htbp]
  \centering
  \begin{subfigure}[b]{0.45\linewidth}
\includegraphics[width=\linewidth]{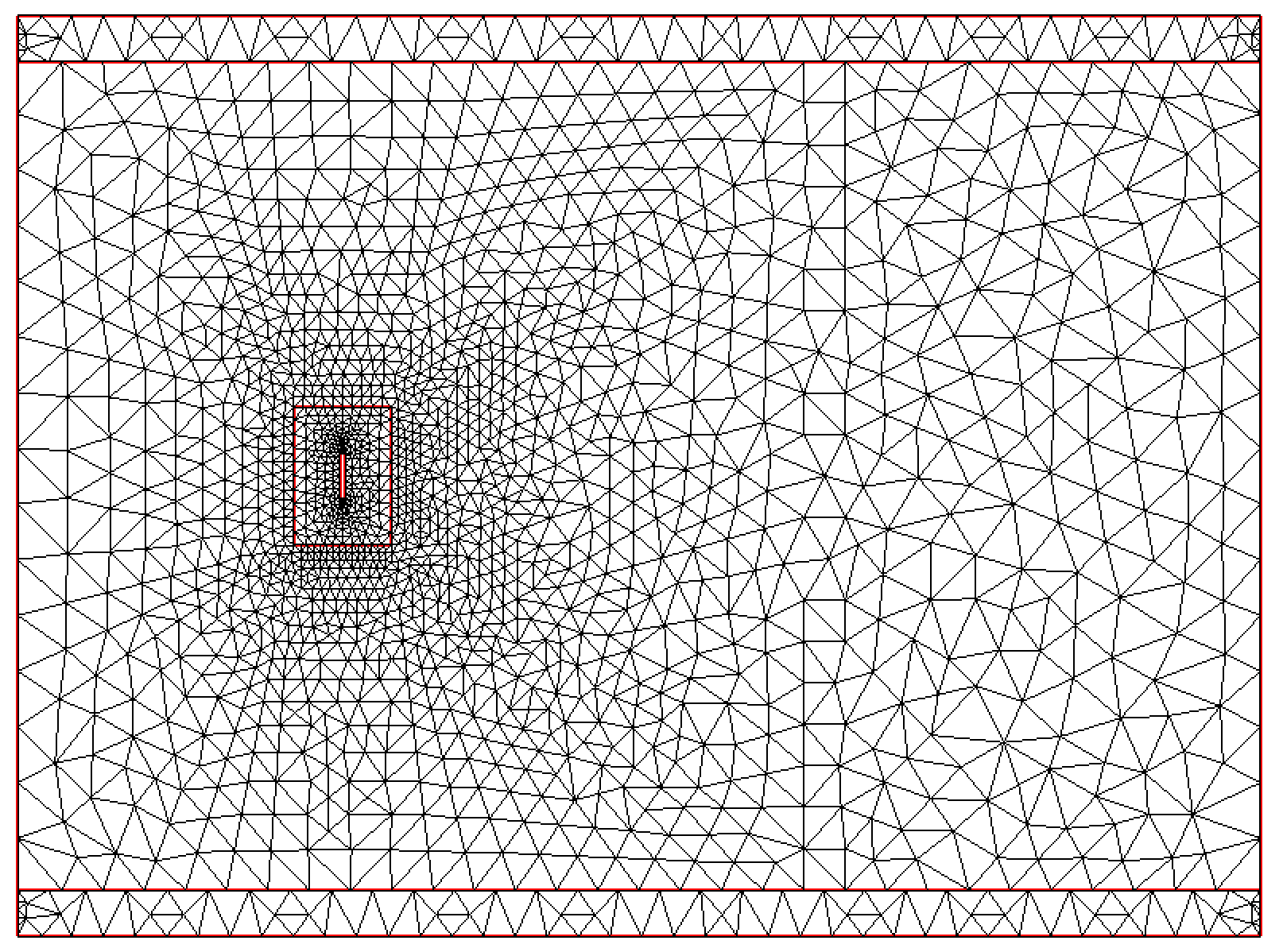}
    \caption{$l=1$.}
    \label{flow_past_obst_mesh-l=1}
  \end{subfigure}
  \begin{subfigure}[b]{0.45\linewidth}
    \includegraphics[width=\linewidth]{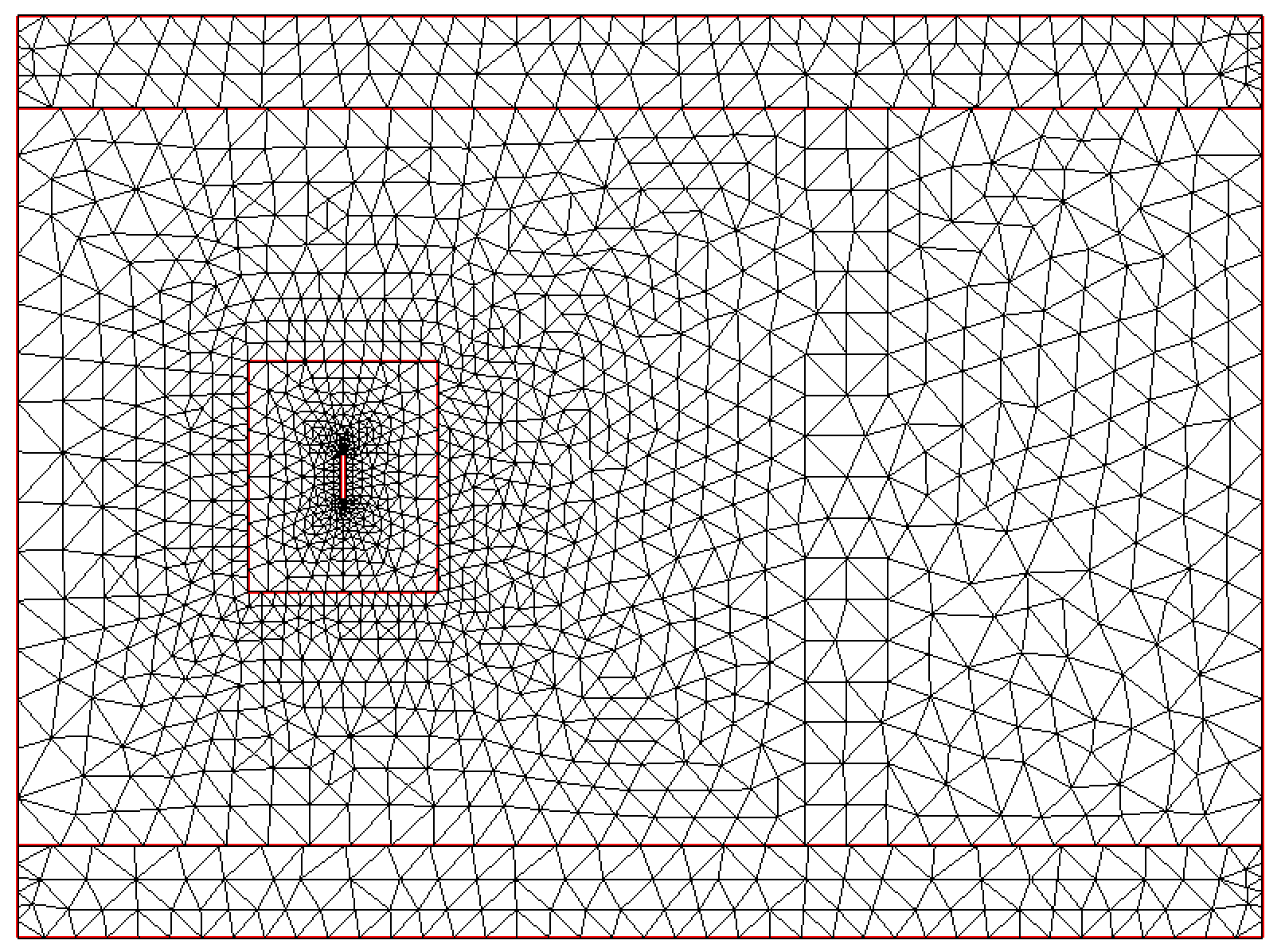}
    \caption{$l=2$.}
     \label{flow_past_obst_mesh-l=2}
  \end{subfigure}
  \begin{subfigure}[b]{0.45\linewidth}
    \includegraphics[width=\linewidth]{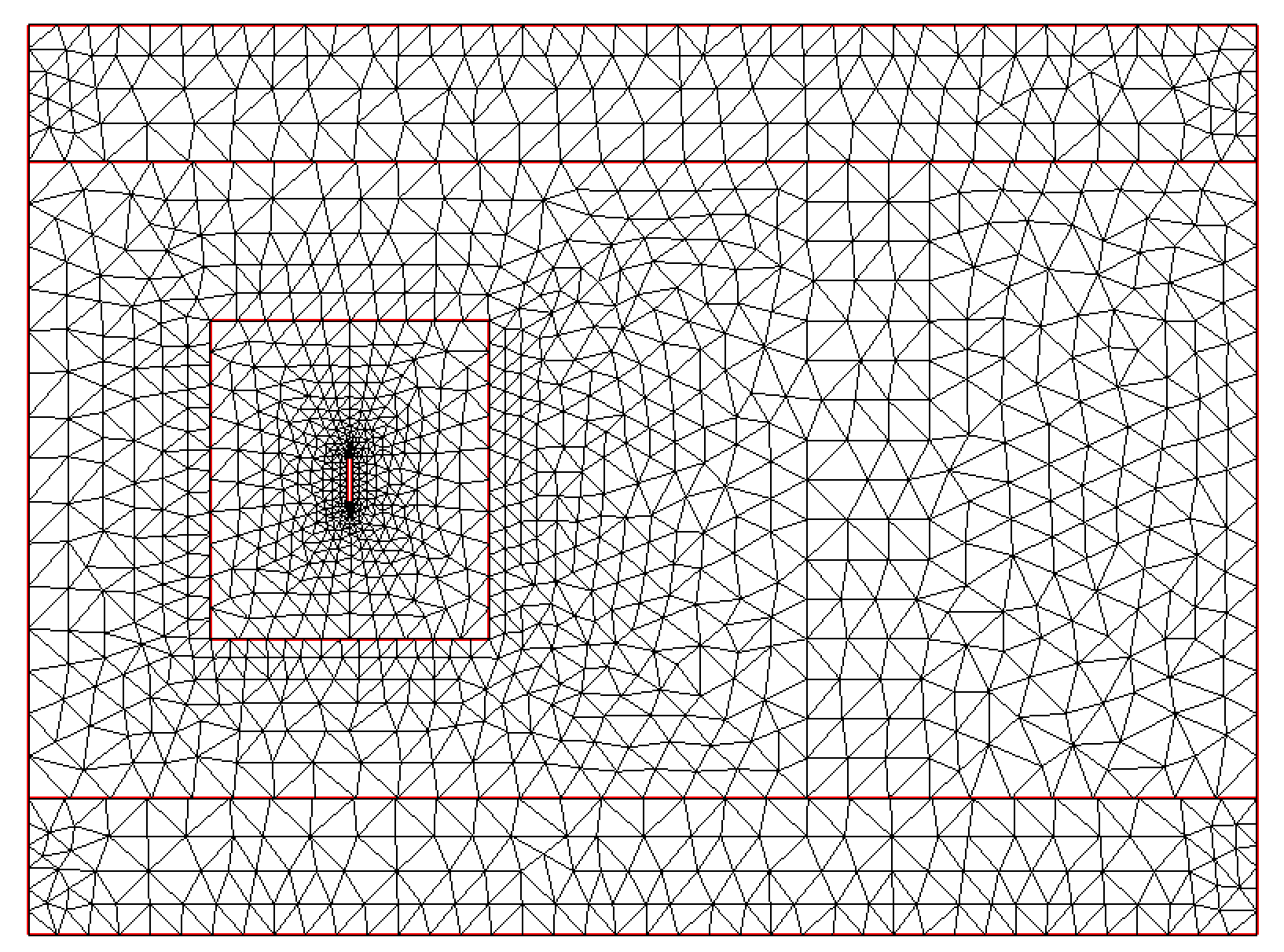}
    \caption{$l=3$.}
     \label{flow_past_obst_mesh-l=3}
  \end{subfigure}
  \caption{Computational Study \RN{1}. Local coarse mesh construction within the domain \(\Omega\), showing three sub-domains.    \(\Omega_{\text{obst}}\):  a box-shaped region centered around the flat plate, 
\(\Omega_{\text{bulk}}\):  the interior region between the top and bottom wall layers and the box-shaped region, \(\Omega_{\text{wall}}\): the boundary layers near the top and bottom walls.
}
\label{flow_past_obst_mesh}
\end{figure}

With the above setup, we examine three cases. First, for thickness parameter \( l = 1 \), we generate a coarse mesh of resolution \( H \) on each subdomain using Delaunay–Voronoi triangulation methods (see Figure~\ref{flow_past_obst_mesh-l=1}), where the observational data are collected. Each triangle in this coarse mesh is then subdivided to create a finer DNS mesh of resolution \( h \); see Table~\ref{table:layer1} for detailed information about this configuration. In addition, to vary the data density and the number of observational data points, we adjust the observed regions \( \Omega_{\text{obst}}, \Omega_{\text{bulk}}, \Omega_{\text{wall}} \) by increasing the layer thickness to \( l = 2 \) and \( l = 3 \) (see Figures~\ref{flow_past_obst_mesh-l=2} and~\ref{flow_past_obst_mesh-l=3}). In each configuration, Regions \(\Omega_{\text{obst}}\) and \(\Omega_{\text{wall}}\) expand, while \(\Omega_{\text{bulk}}\) shrinks accordingly. See Tables~\ref{table:layer2} and~\ref{table:layer3} for further details.

\begin{table}[ht!]
\centering
\begin{tabular}{|c c c c|} 
 \hline
 Domain & $\#$ Vertices &    $\#$ Vertices per Area  & Max size of mesh   \\ [0.5ex] 
 \hline
$\Omega$ (DNS) & 7322  & 14 & 0.85 \\ 
$\Omega_{\text{obst}}$ & 278 & 45  & 0.40\\
$\Omega_{\text{bulk}}$ & 1447  & 3  &  1.70 \\
$\Omega_{\text{wall}}$ & 204 & 4  &  1.16  \\ [1ex] 
 \hline
\end{tabular}
\caption{Computational Study \RN{1}. Details on the flow over a flat obstacle meshes for $l=1$.}
\label{table:layer1}
\end{table}

\begin{table}[ht!]
\centering
\begin{tabular}{|c c c c|} 
 \hline
 Domain & $\#$ Vertices &    $\#$ Vertices per Area  & Max size of mesh   \\ [0.5ex] 
 \hline
$\Omega$ (DNS) & 6312  & 12 & 0.75 \\ 
$\Omega_{\text{obst}}$ & 364 & 18  & 0.68\\
$\Omega_{\text{bulk}}$ & 1003  & 2 &  1.51 \\
$\Omega_{\text{wall}}$ & 310 & 3  &  1.22 \\ [1ex] 
 \hline
\end{tabular}
\caption{Computational Study \RN{1}. Details on the flow over a flat obstacle meshes for $l=2$.}
\label{table:layer2}
\end{table}

\begin{table}[ht!]
\centering
\begin{tabular}{|c c c c|} 
 \hline
 Domain & $\#$ Vertices &    $\#$ Vertices per Area  & Max size of mesh   \\ [0.5ex] 
 \hline
$\Omega$ (DNS) & 5736  & 11 & 0.69\\ 
$\Omega_{\text{obst}}$ & 428 & 10 & 0.95\\
$\Omega_{\text{bulk}}$ & 722 & 2  &  1.38 \\
$\Omega_{\text{wall}}$ & 383 & 2  &  1.35  \\ [1ex] 
 \hline
\end{tabular}
\caption{Computational Study \RN{1}. Details on the flow over a flat obstacle meshes for $l=3$.}
\label{table:layer3}
\end{table}
We do not impose an external force, and the flow is driven by the inflow velocity, $u_{\text{in}} = \langle 1, 0 \rangle^T$. We impose the no-slip boundary conditions at the obstacle and walls, and a weak zero-traction boundary condition at the outflow. We set the Reynolds number $\Rey \sim  600$, a time step of $\Delta t = 0.02$, and a final time of $T = 81$, corresponding to three turnover times. Here, the turnover time is based on the full streamwise domain length, long wall length $= 27$. Since this test examines a through-flow problem, it evaluates errors over three turnover times rather than an extended period. Nevertheless, the flow is complicated, with many intriguing features, and is a standard test problem \cites{GLRVZ21, CIBIK2025117526}. 

The true initial condition $\bfu_0$ is generated by solving the steady Stokes problem with the same boundary conditions and zero external force, and then the true solution is generated by solving the NSE on the DNS global mesh h. The total degrees of freedom for the global fine mesh $50182$.  To initiate the data assimilation simulation, we assume that the flow in the DA system is initially at rest, i.e., \( \bfv_0 = 0 \) at \( t = 0 \), and we set the nudging parameter to \( \mu = 5 \).

\subsection*{Impact of Observation Placement: Boundary vs. Interior}

\begin{figure}[ht!]
  \centering
  \begin{subfigure}[b]{0.49\linewidth}
    \includegraphics[width=0.8\linewidth]{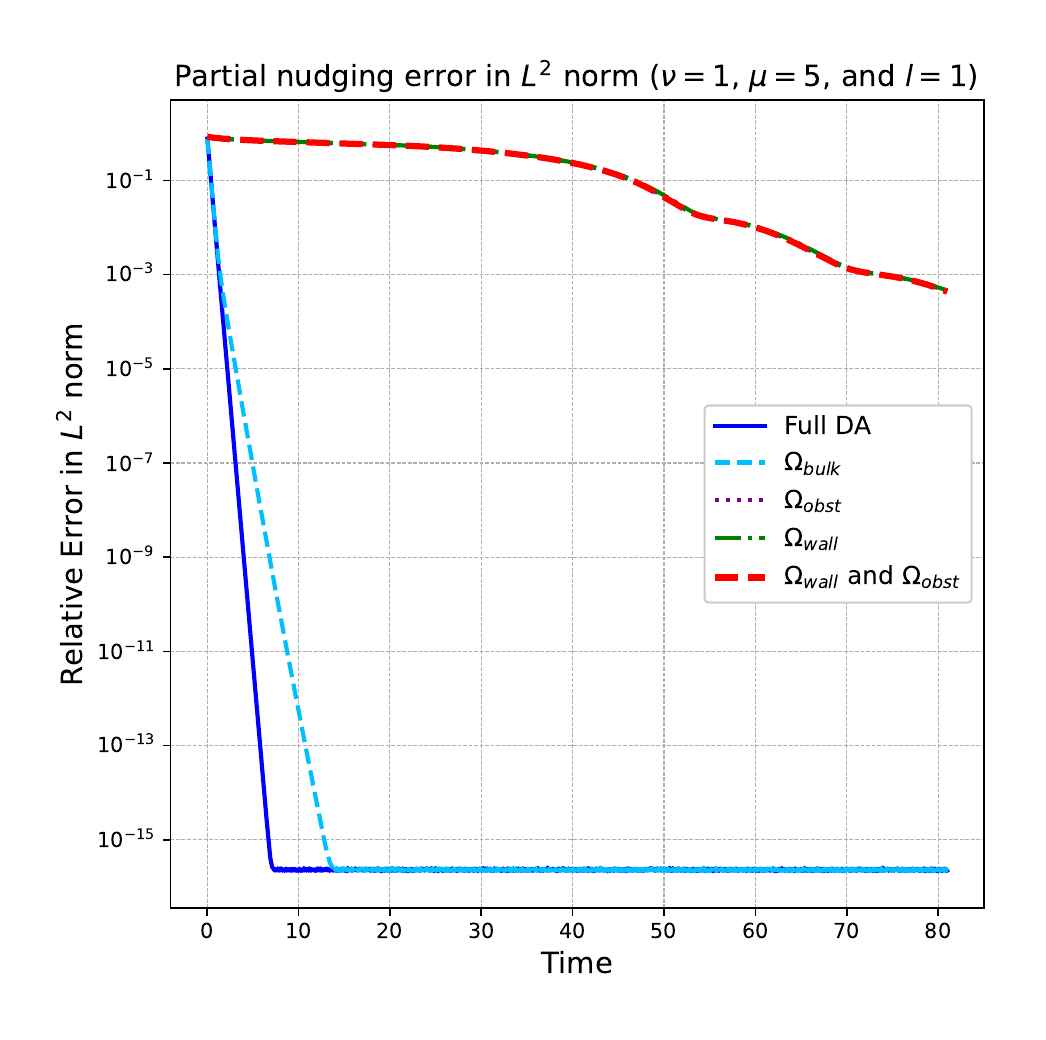}
     \caption{Flow over a flat obstacle with $l=1$.}
     \label{fig:layer1_L2_error}
  \end{subfigure}
  \begin{subfigure}[b]{0.49\linewidth}
    \includegraphics[width=0.8\linewidth]{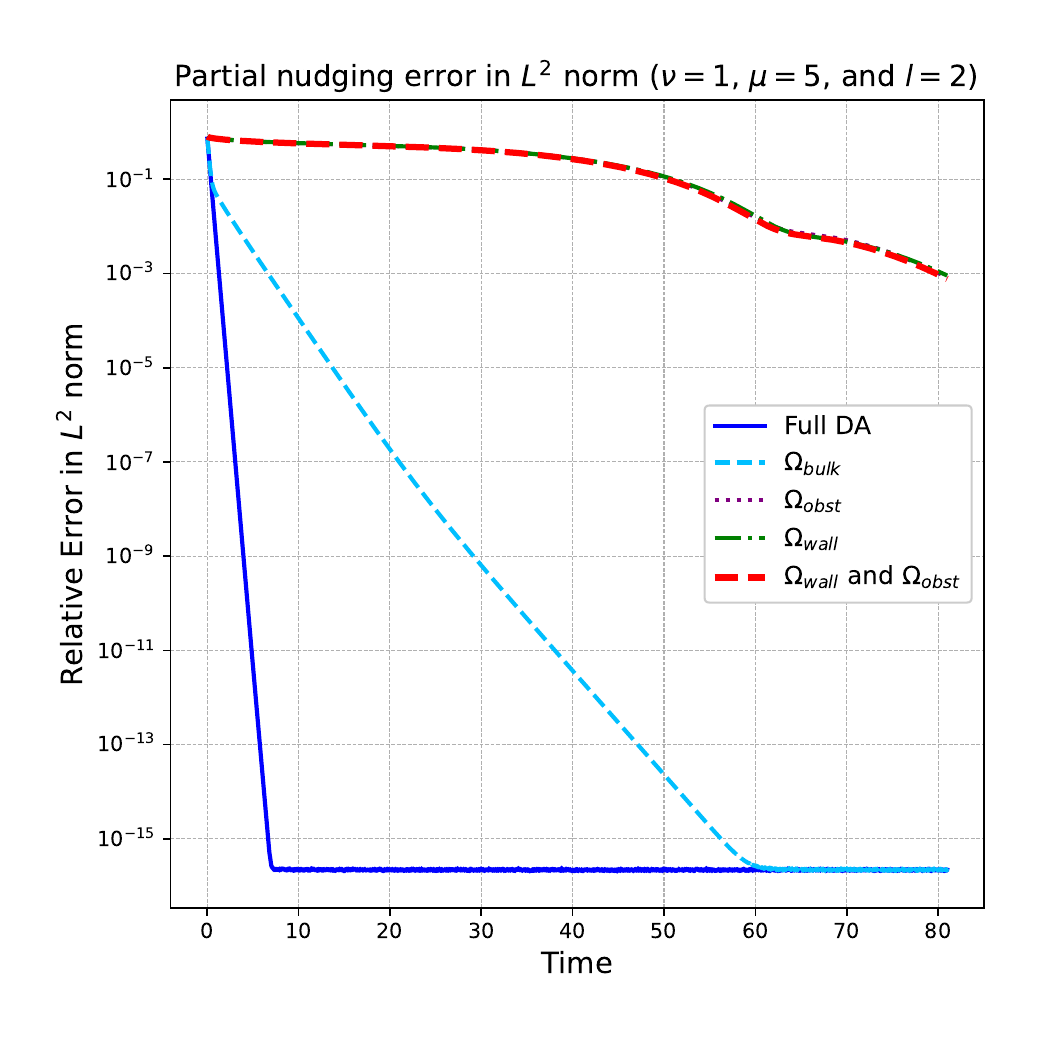}
 \caption{Flow over a flat obstacle with $l=2$.}
 \label{fig:layer2_L2_error}
  \end{subfigure}
   \begin{subfigure}[b]{0.49\linewidth}
    \includegraphics[width=0.8\linewidth]{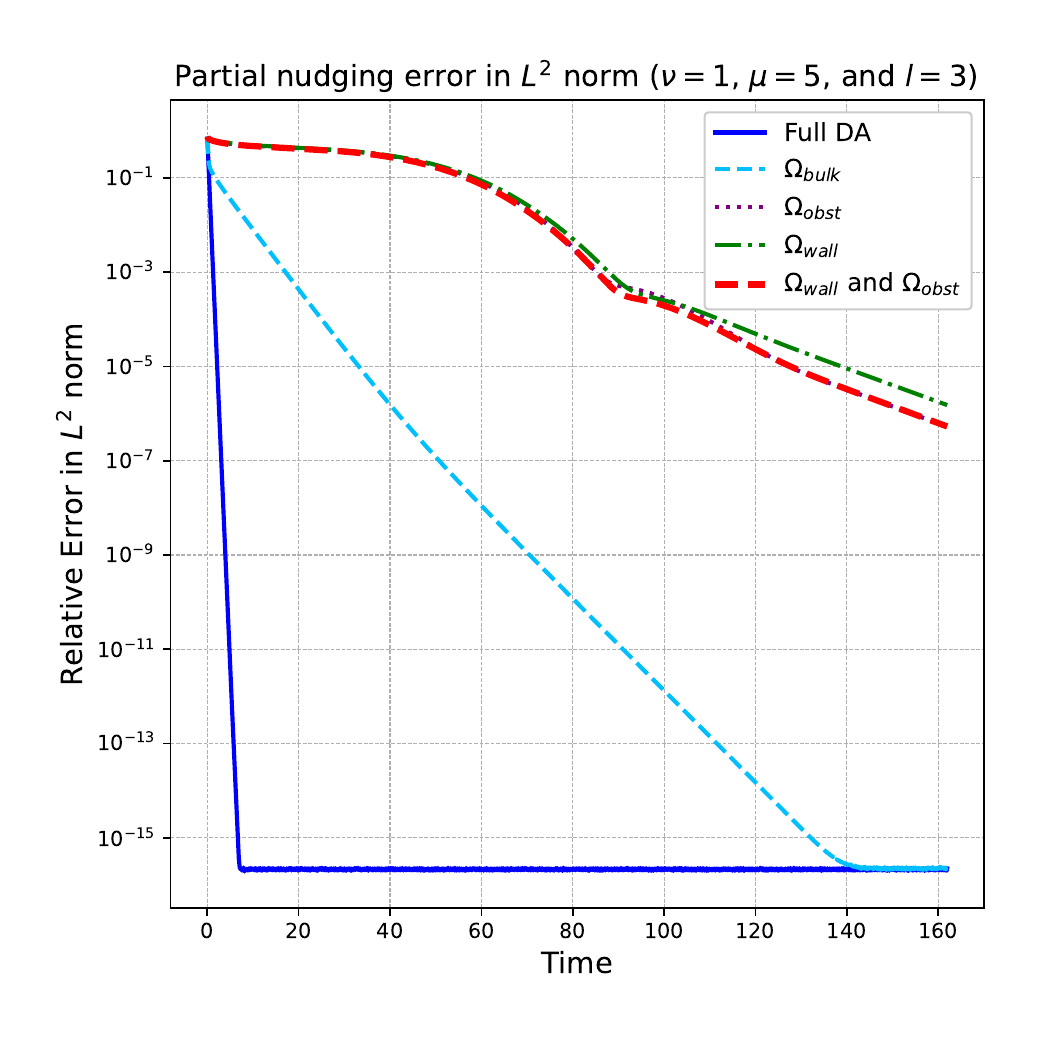}
    \caption{Flow over a flat obstacle with $l=3$.}
    \label{fig:layer3_L2_error}
  \end{subfigure}
  \caption{Computational Study \RN{1}. Comparison of results using the local  data assimilation algorithm applied to the flow over a flat obstacle with data collected from different regions of the domain: $\Omega_{\text{obst}}$, $\Omega_{\text{bulk}}$, and $\Omega_{\text{wall}}$.}
   \label{fig:L2_error_Obstacle}
  \end{figure}

Figure~\ref{fig:L2_error_Obstacle} shows the time evolution of the relative error in the \( L^2 \) norm of the velocity difference between the reference solution \( \bfu \) and the locally nudged solutions \( \bfv \), with observational data collected from different regions:  Obstacle-adjacent Region (\( \Omega_{\text{obst}} \)), Interior bulk Region (\( \Omega_{\text{bulk}} \)), Wall-adjacent Region (\( \Omega_{\text{wall}} \)), and the union of \(\Omega_{\text{obst}} \cup \Omega_{\text{wall}} \). 

For comparison, the results of global (full-domain) nudging are also shown in each case (solid blue curve), where data are collected throughout the entire domain \( \Omega \). Across all three cases (\( l = 1, 2, 3 \)), full synchronization is achieved when using data from the interior region \( \Omega_{\text{bulk}} \), which is located away from the boundaries, though the convergence occurs slightly more slowly than in the global data assimilation case.  In contrast, data collected near the boundaries (i.e., from Regions \( \Omega_{\text{obst}} \) and \( \Omega_{\text{wall}} \)) contribute little to no improvement in synchronization. This limited effectiveness persists even when observations from \(\Omega_{\text{obst}} \cup \Omega_{\text{wall}} \) are combined.  Note that in our final simulation (Figure~\ref{fig:layer3_L2_error}), where the interior observation region $\Omega_{\text{bulk}}$ was reduced, global synchronization was still achieved. However, the convergence required considerably more time, underscoring that while interior observations are indeed sufficient, the spatial extent of the observed region plays a pivotal role in determining the efficiency of the assimilation process.

\begin{figure}[htbp]
    \centering
\includegraphics[width=0.45\linewidth]{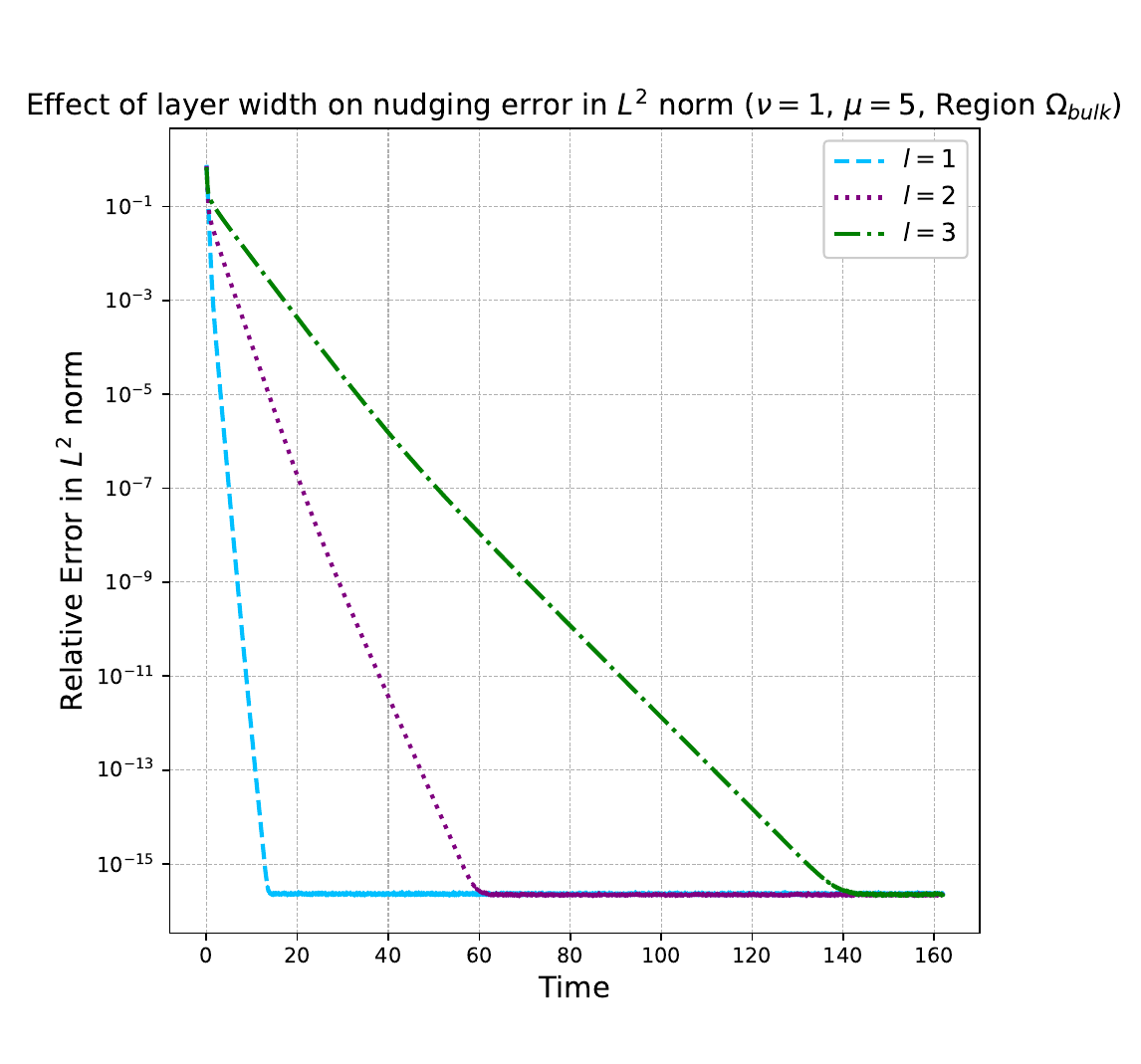}  
\caption{Computational Study \RN{1}. Relative error when applying local DA to $\Omega_{\text{bulk}}$ with different $l$.}
\label{fig:layers_L2}
\end{figure}

In addition, in Figure~\ref{fig:layers_L2}, we compare the performance of the local nudging algorithm focusing on $\Omega_{\text{bulk}}$, for different values of the layer thickness parameter \( l \), over a longer time interval \( T = 182 \). The error corresponding to \( l = 1 \) decreases the fastest, reaching a full synchronization around \( t = 15 \). For \( l = 2 \), convergence occurs around \( t = 60 \), and for \( l = 3 \), the full synchronization is achieved near \( t = 140 \). This trend is due to the fact that expanding the observational region provides more information, which improves the assimilation performance.

\subsection{Local DA with different nudging parameter}
\begin{figure}[htbp]
    \centering
\includegraphics[width=0.45\linewidth]{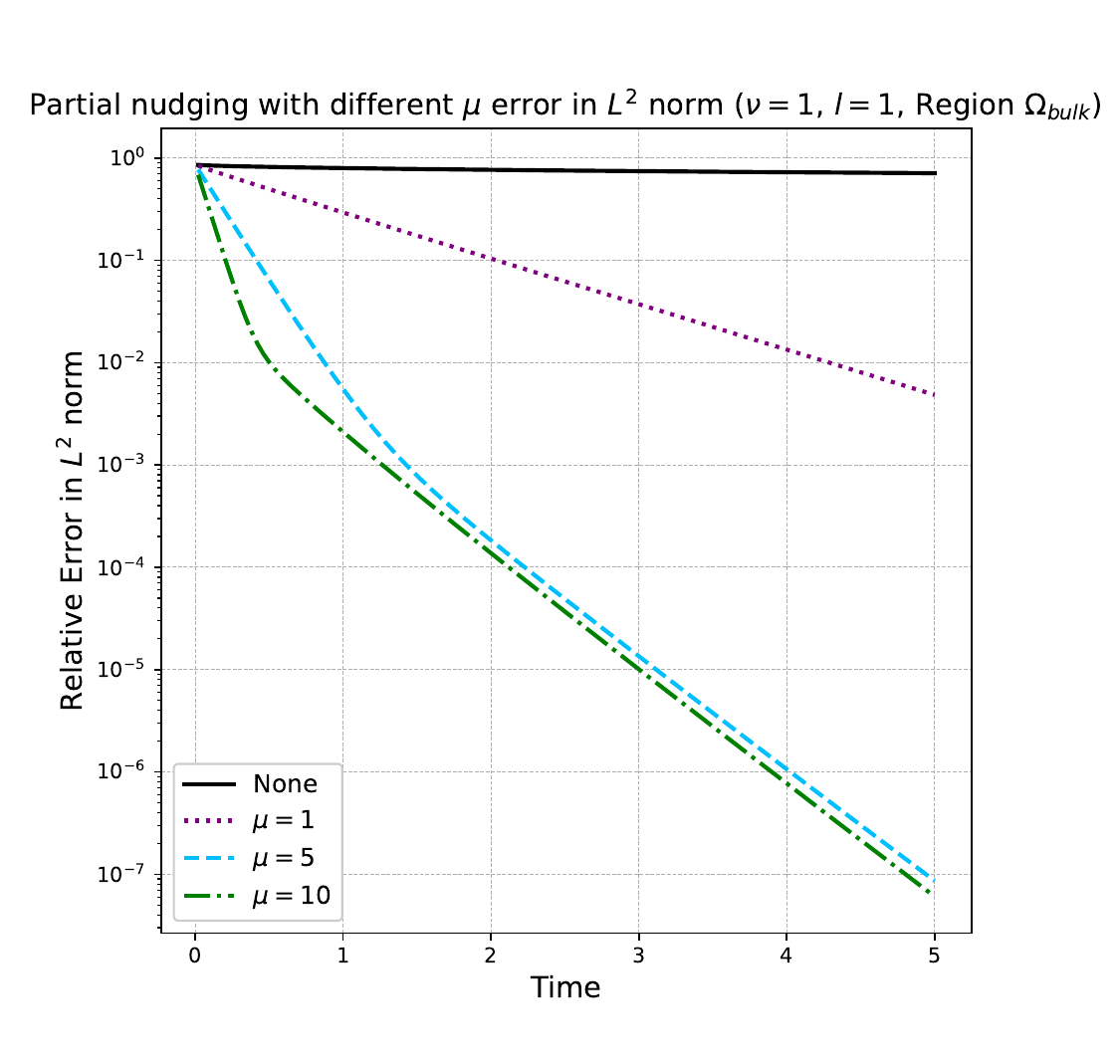} 
\caption{Computational Study \RN{1}. Relative error when applying DA to $\Omega$ with different $\mu$ values.}
\label{fig:Regions_L2_mu}
\end{figure}

We also vary the nudging parameter $\mu$ applied to Interior bulk Region, with $\mu = 0, 1, 5, 10$ to study the its effect on the error. We set $l = 1$ for all experiments. In Figures \ref{fig:Regions_L2_mu}, we observe that DA applied to Interior bulk Region shows significantly better results for larger $\mu$, with $\mu = 5$ and $\mu = 10$ achieving similar results, over a longer period. For $\mu=1$, a smaller value, the performance is much better than the no-DA system.

\subsection{Local DA with noisy data}

In this subsection, we study the effects of noise in the interior observational
data (DA interior bulk region) on local nudging. The noise is added to the
noise-free DNS solution. We denote the noise amplitude by
\(\delta_{\mathrm{noise}}\). Specifically, the noise is sampled from a uniform
distribution on $[-1,1]$ and scaled by \(\delta_{\mathrm{noise}}\). We use the
same numerical setting as in the previous subsection and consider \(l=1\),
\(\nu=1\), and \(\mu=5\).
We compare local nudging with noise-free interior data, \(\delta_{\mathrm{noise}}=0\), and noisy interior data with noise amplitudes \(\delta_{\mathrm{noise}}=0.001\), \(0.01\), and \(0.05\).

\begin{figure}[ht!]
    \centering
    \includegraphics[width=0.43\linewidth]{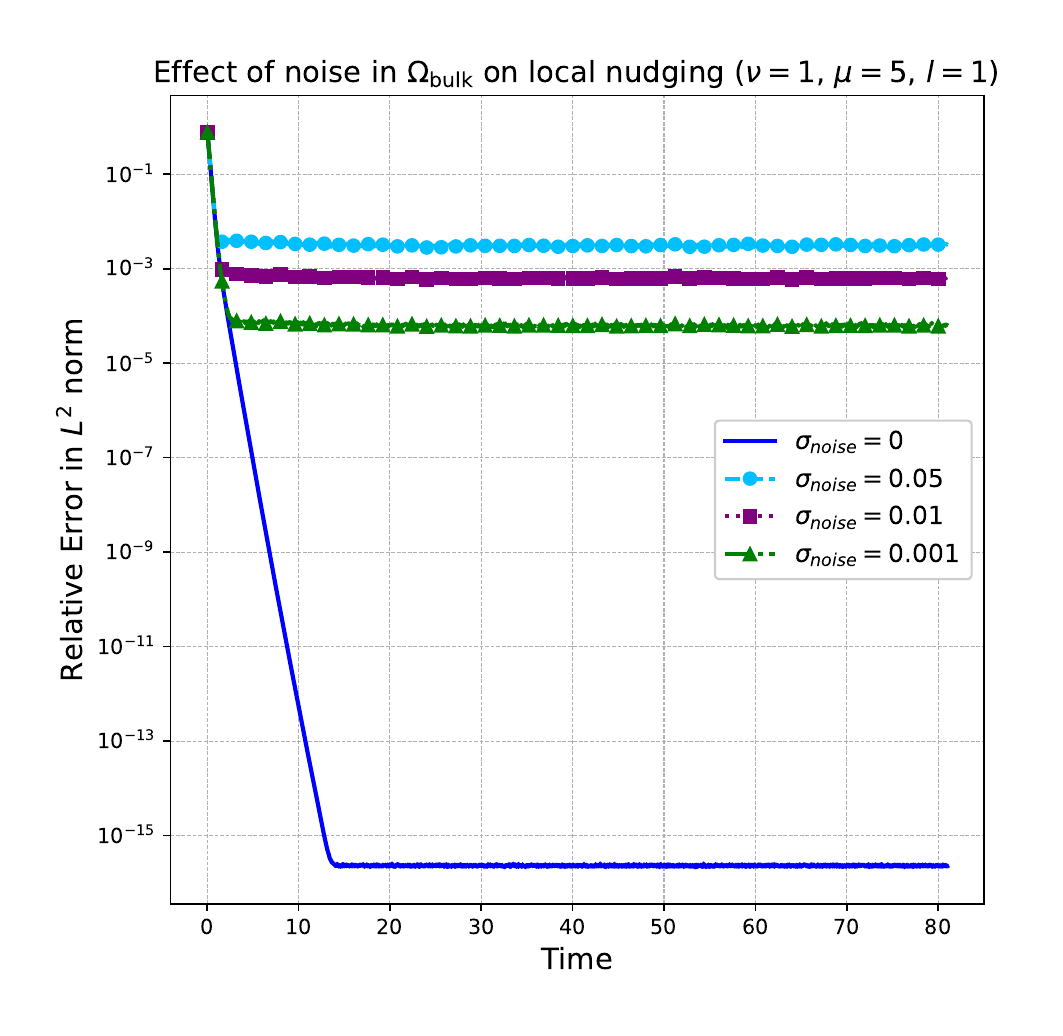}
    \caption{Relative error plot of velocity in $L^2$ norm with different noise amplitudes $\delta_{\mathrm{noise}}=0.05, 0.01,0.001$ and $0$.}
    \label{fig:noise_level}
\end{figure}

Figure \ref{fig:noise_level} shows the relative velocity error in the \(L^2\) norm for different noise amplitudes, \(\delta_{\mathrm{noise}}= 0, 0.001, 0.01,\) and \(0.05\). With noise-free interior observational data, the numerical solution achieves full synchronization, with the relative error decreasing to approximately \(10^{-16}\). In the presence of observational noise, the error initially decreases but eventually reaches a noise-dependent plateau. In particular, the relative error is approximately \(10^{-5}\), \(10^{-4}\), and \(10^{-3}\) for \(\delta_{\mathrm{noise}}=0.001\), \(0.01\), and \(0.05\), respectively. Thus, increasing the noise amplitude leads to a correspondingly higher asymptotic error level.

\begin{figure}[ht!]
    \centering
    \includegraphics[width=0.45\linewidth]{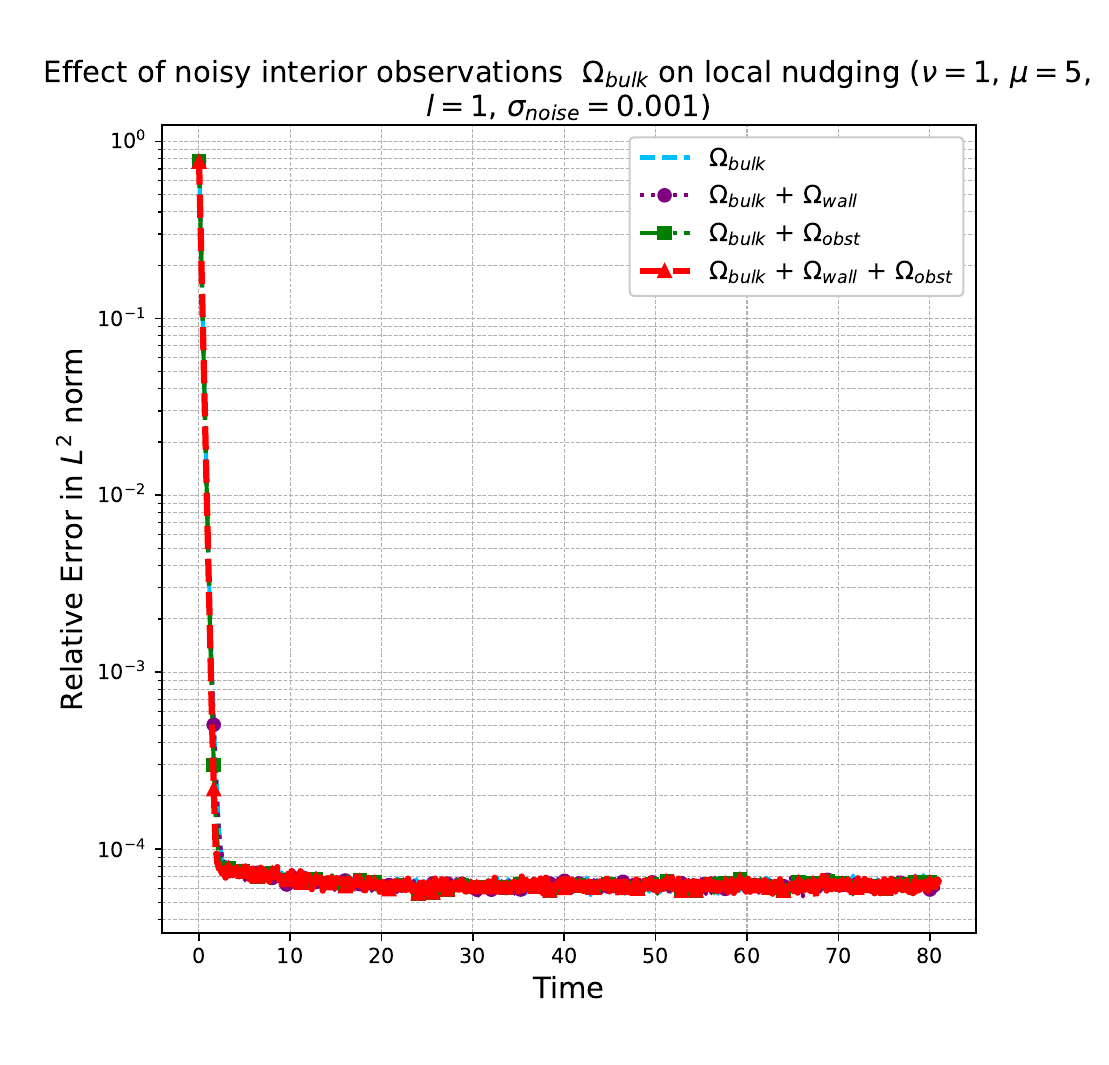}
    \caption{Relative error plot of velocity in $L^2$ norm with noise amplitudes $\delta_{\mathrm{noise}}=0.001$ for interior data combined with boundary data.}
    \label{fig:noise_regions}
\end{figure}

We investigate whether adding boundary observations to the noisy interior data helps reduce the assimilation error. We compare noisy interior observations alone with three additional observation configurations: interior with wall observations, interior with obstacle-boundary observations, and interior with both wall and obstacle-boundary observations. We set \(\delta_{\mathrm{noise}}=0.001\), \(\nu=1\), \(\mu=5\), and \(l=1\). Figure \ref{fig:noise_regions} shows that adding boundary observations does not significantly reduce the assimilation error in the presence of noise in the interior observational data.

\section{Computational Study II: The Couette flow}\label{Section; SIM2}
In this section, we examine two-dimensional  Couette flow with tangential boundary conditions in the absence of a body force. Our numerical experiment is inspired by the work of Couette and Taylor, known as the Taylor experiment \cite{F95}. In Taylor’s experiments, two concentric, very long cylinders are used, with the space between them filled with a liquid. Each cylinder can rotate independently. Consider a $2D$ slice of this experiment, where the domain \(\Omega\) is an annulus defined as
$$
\Omega = \{ (x,y) \in \mathbb{R}^2 : 0.1^2 \leq x^2 + y^2 \leq 1^2 \}.
$$
The flow is driven by rotating both the inner and outer circles, with one rotating clockwise and the other counterclockwise. In the absence of a body force $\bff= \bfzero$, we consider the Reynolds number to be $\Rey= 600$.

\begin{figure}[ht!]
  \centering
  \begin{subfigure}[b]{0.49\linewidth}
    \includegraphics[width=0.8\linewidth]{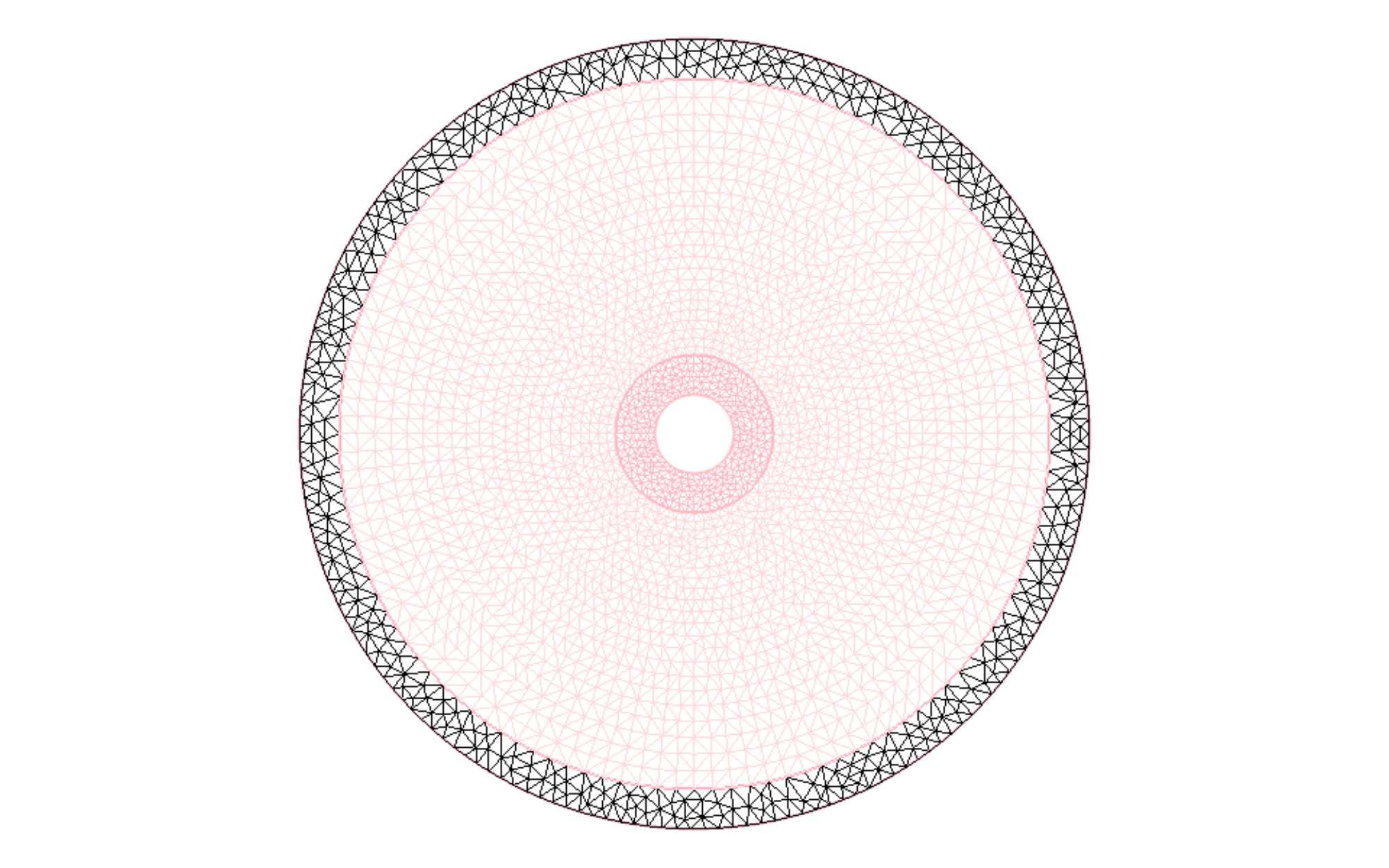}
    \caption{$r_1=0.2, \ r_2=0.9$.}
    \label{fig: sec6_r1_point2_r2_point9}
  \end{subfigure}
  \begin{subfigure}[b]{0.49\linewidth}
    \includegraphics[width=0.8\linewidth]{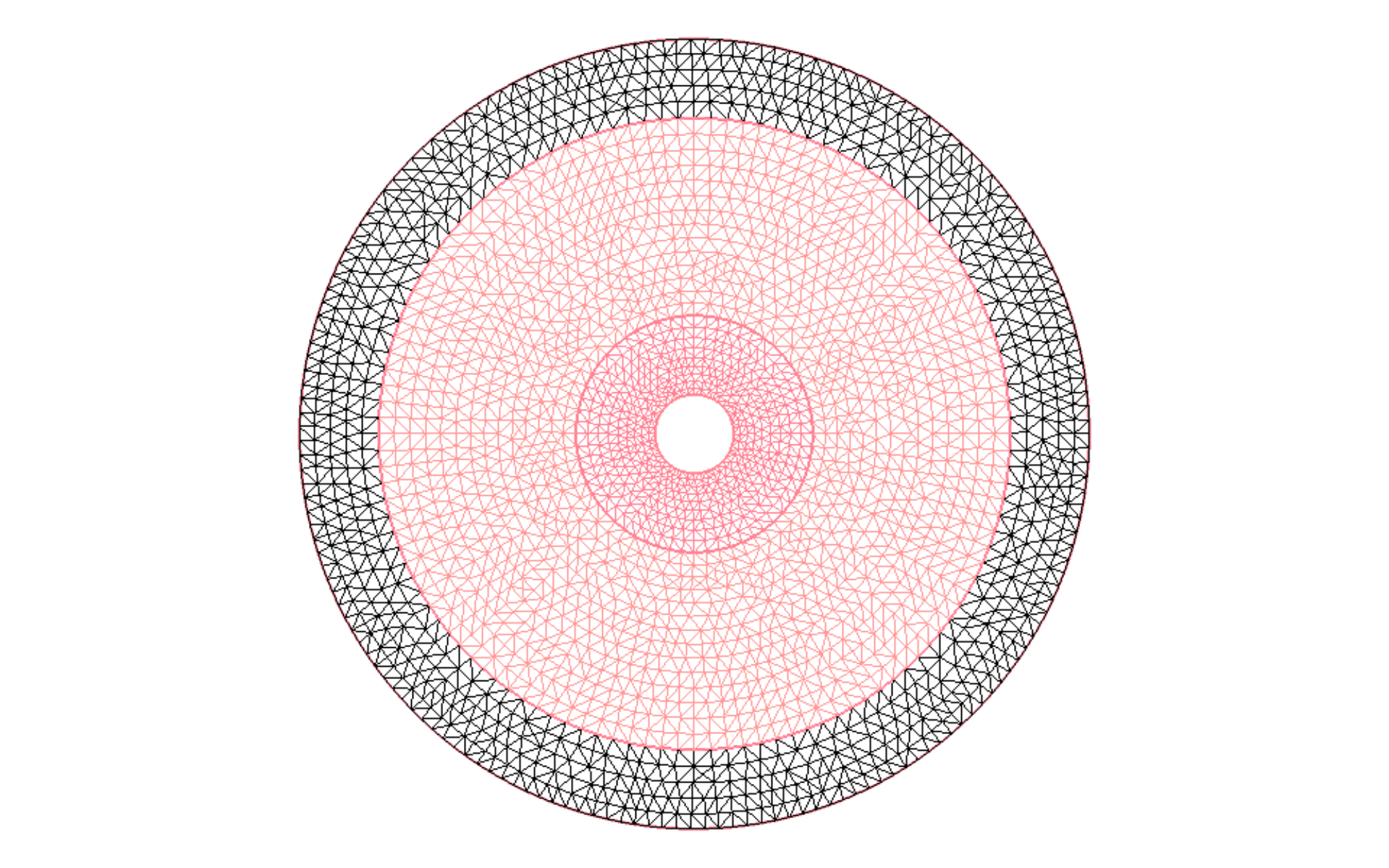}
    \caption{$r_1=0.3, \ r_2=0.8$.}
    \label{fig: r1_point3_r2_point8}
  \end{subfigure}
   \begin{subfigure}[b]{0.49\linewidth}
    \includegraphics[width=0.8\linewidth]{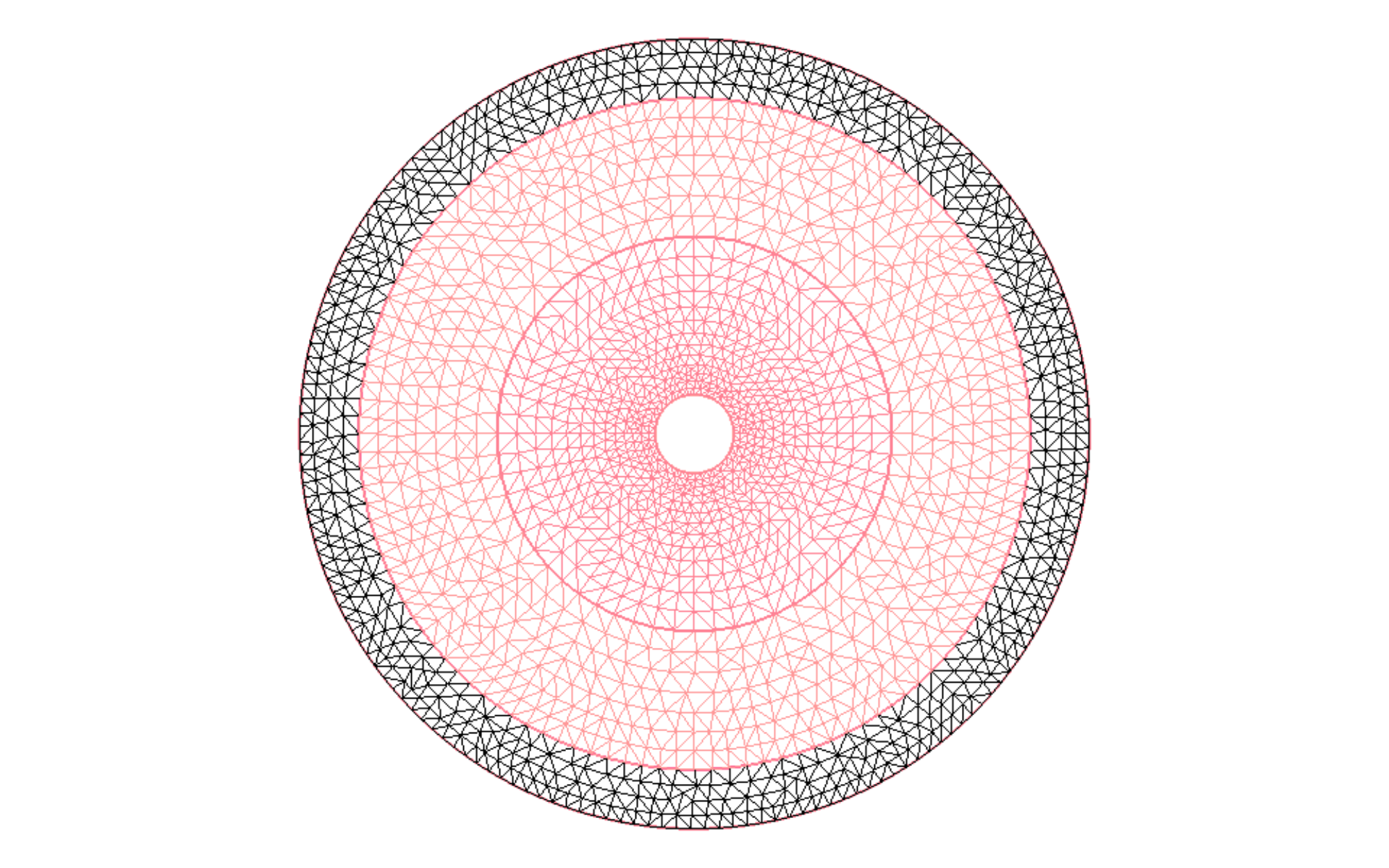}
    \caption{$r_1=0.5, \ r_2=0.85$.}
    \label{fig: r1_point5_r2_point85}
  \end{subfigure}
  \caption{Computational Study \RN{2}. Local coarse mesh construction within the domain \(\Omega\) showing three sub-domains: \(\Omega_{\text{obst}}\) (closest to the inner disc - dark red), \(\Omega_\text{bulk}\) (between the inner and outer regions - light red), and \(\Omega_\text{wall}\) (closest to the outer boundary - black).}
  \label{fig:LocMeshShear}
\end{figure}

Within the domain $\Omega$, we introduce two artificial circular boundaries
with radii
\[
r_1 = 0.2
\qquad\text{and}\qquad
r_2 = 0.9.
\]
For $s\in[0,2\pi]$, the physical and artificial boundaries are parametrized by
\begin{align*}
C_1(s) &= \bigl(0.1\cos(s),\,0.1\sin(s)\bigr), \\
C_2(s) &= \bigl(r_1\cos(s),\,r_1\sin(s)\bigr), \\
C_3(s) &= \bigl(r_2\cos(s),\,r_2\sin(s)\bigr), \\
C_4(s) &= \bigl(\cos(s),\,\sin(s)\bigr).
\end{align*}
Here, $C_1$ and $C_4$ represent the inner and outer physical boundaries, respectively, whereas $C_2$ and $C_3$ are artificial boundaries used to partition the domain. Using these boundaries, we decompose $\Omega$ into three physically distinct subdomains:
\[
\Omega
=
\Omega_{\mathrm{obst}}
\cup
\Omega_{\mathrm{bulk}}
\cup
\Omega_{\mathrm{wall}},
\]
whose interiors are pairwise disjoint; see
Figure~\ref{fig:LocMeshShear}. The subdomains are defined as follows:
\begin{itemize}
    \item  $\Omega_{\mathrm{obst}}$: Obstacle-adjacent region, the annular region enclosed between $C_1$ and $C_2$, adjacent to the inner disc.

    \item $\Omega_{\mathrm{bulk}}$: Interior bulk region, the annular region enclosed between $C_2$ and $C_3$, located away from
    both physical boundaries.

    \item $\Omega_{\mathrm{wall}}$: Wall-adjacent region, the annular region enclosed between $C_3$ and $C_4$, adjacent to the outer wall.
\end{itemize}
For convenience, we also define the complete boundary-adjacent region by
\[
\Omega_{\mathrm{bdry}}
:=
\Omega_{\mathrm{obst}}
\cup
\Omega_{\mathrm{wall}}.
\]

We first construct a coarse mesh $H$ on each sub-domain using Delaunay-Voronoi methods, see Figure \ref{fig: sec6_r1_point2_r2_point9}, where the data are collected. Each triangular mesh is then split into four sub-meshes to generate the DNS mesh $h$. Throughout these  refinements,  we ensure that scales are resolved down to the Kolmogorov dissipation scale, $h_{\text{max}} \leq \Rey^{-\frac{1}{2}}$. Detailed information about both the DNS mesh and the local meshes can be found in Table \ref{table:1}. 

\begin{table}[ht!]
\centering
\begin{tabular}{|c c c c|} 
 \hline
 Domain & $\#$ Vertices &    $\#$ Vertices per Area  & Max size of mesh   \\ [0.5ex] 
 \hline
$\Omega$ (DNS) & 12293  & 3965 & 0.037   \\ 
$\Omega_{\text{obst}}$ & 356 & 3747  & 0.029 \\
$\Omega_{\text{bulk}}$ & 2379  & 992  &  0.064 \\
$\Omega_{\text{wall}}$ & 572 & 969  &  0.074  \\ [1ex] 
 \hline
\end{tabular}
\caption{Computational Study \RN{2}. Details on the Couette spatial meshes for $r_1=0.2, \ r_2=0.9$.}
\label{table:1}
\end{table}

Using the schemes described in Algorithm \ref{Algorithm}, we consider the time interval $[0, 100]$ with a time-step size $\Delta t = 0.01$, following the approach in \cite{JL14}. To generate the initial condition for the true solution, we begin by running the simulation with the specified shear boundary condition from an initial state of zero up to \( t = 5 \). The resulting state at \( t = 5 \) is then used as the initial condition for subsequent simulations, denoted as \(\bfu^h_0\). For the DA algorithms, we start from $\bfv^h_0=\bfzero$ with $\mu=10$.  For the interpolant, we use local linear interpolation on a mesh \( H \) within the local region $\Omega_{\text{bulk}}$, which is four levels coarser than the DNS mesh.  

\subsection*{Impact of Observation Placement: Boundary vs. Interior}

\begin{figure}[htbp]
  \centering
  \begin{subfigure}[b]{0.45\linewidth}
\includegraphics[width=0.75\linewidth]{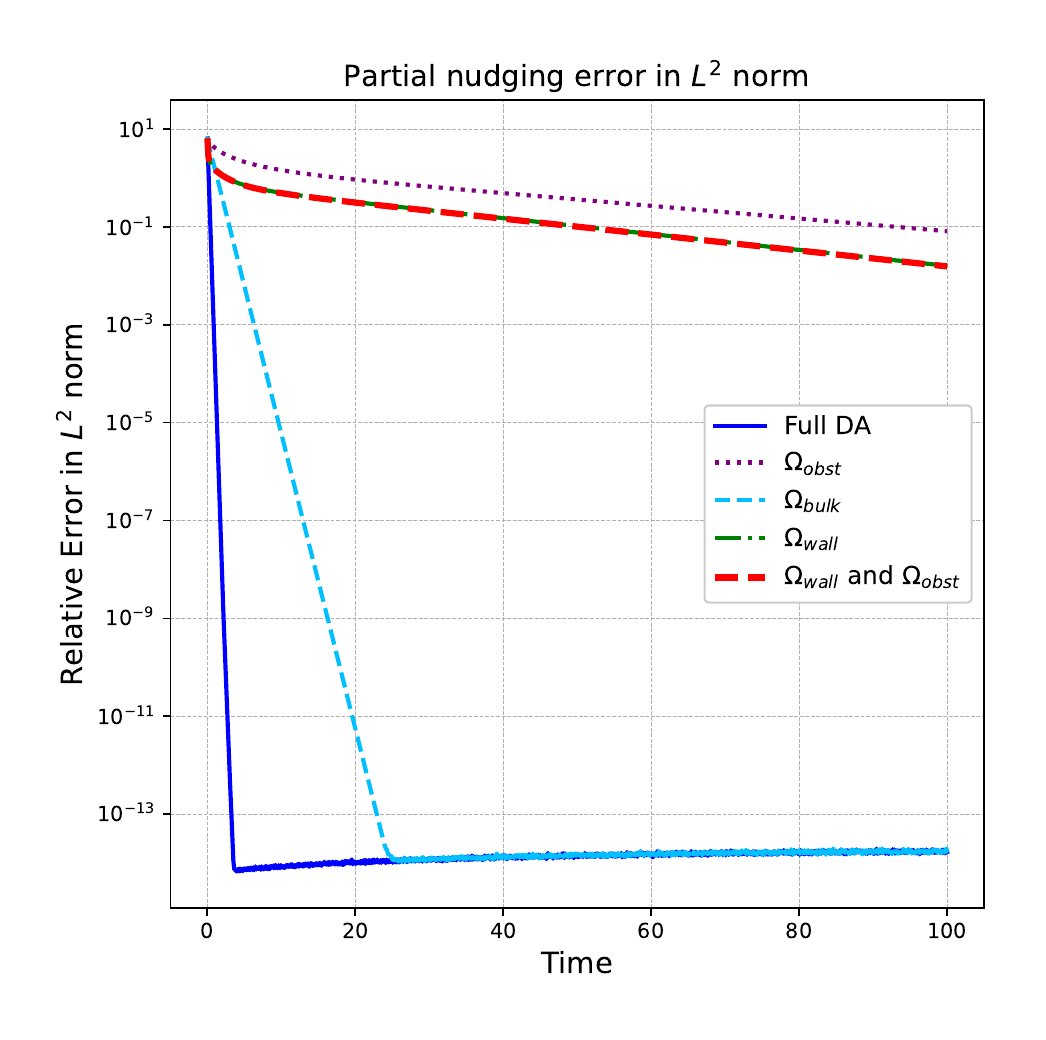}
    \caption{Couette flow $r_1=0.2,\  r_2=0.9$.}
    \label{fig: r1_point2_r2_point9}
  \end{subfigure}
  \begin{subfigure}[b]{0.45\linewidth}
\includegraphics[width=0.75\linewidth]{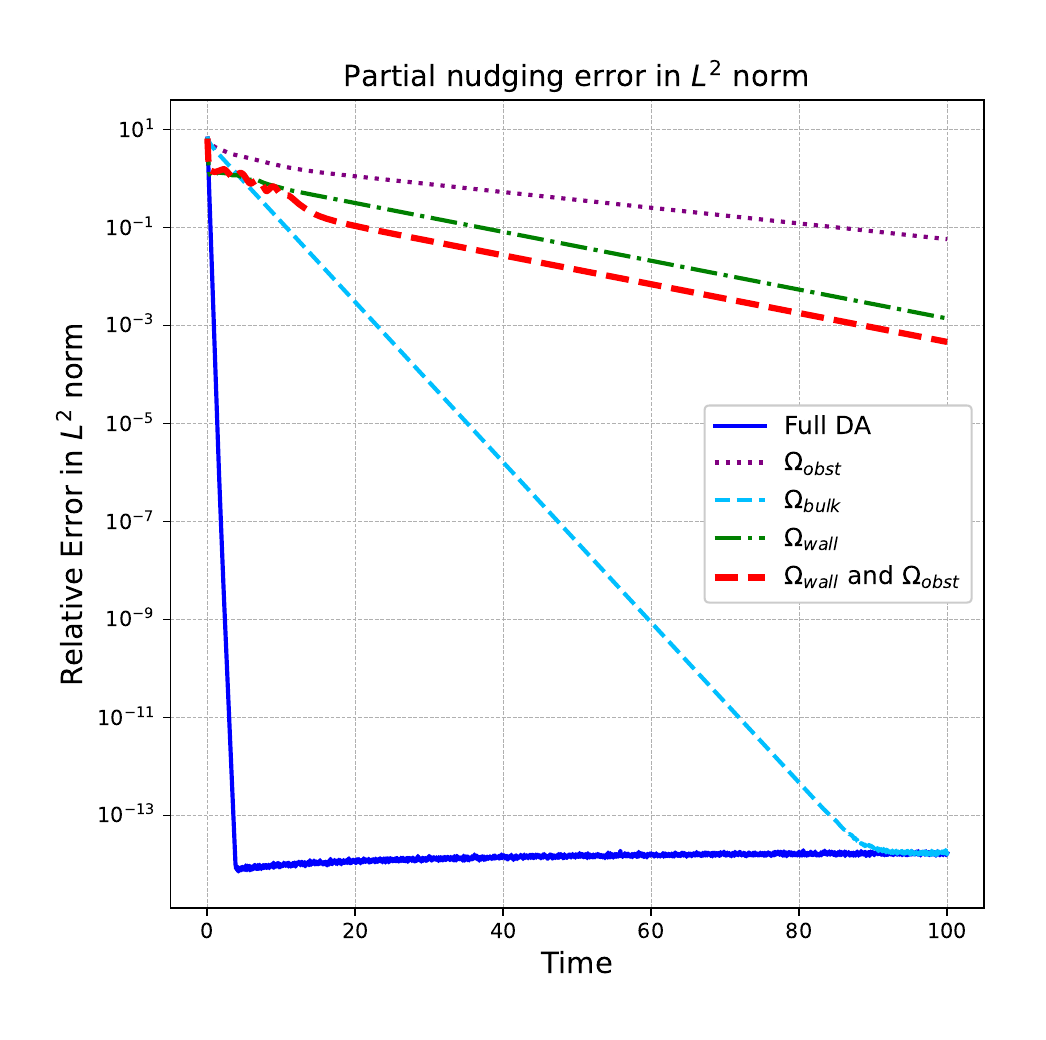}
    \caption{Couette flow $r_1=0.3,\  r_2=0.8$.}
    \label{fig: r1_point2_r2_point10}
  \end{subfigure}
   \begin{subfigure}[b]{0.45\linewidth}    \includegraphics[width=0.75\linewidth]{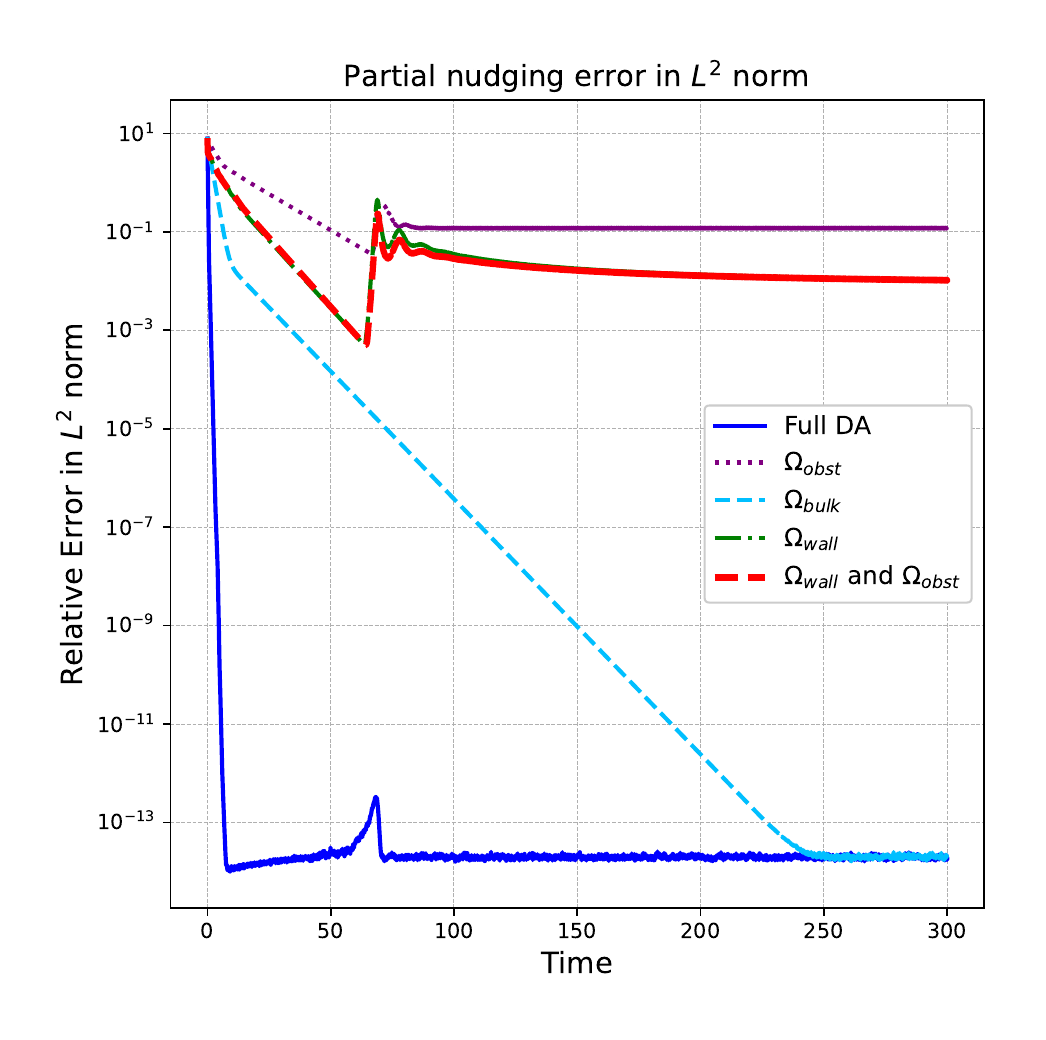}
    \caption{Couette flow $r_1=0.5, \  r_2=0.85$.}
    \label{fig: r1_point2_r2_point11}
  \end{subfigure} 
\caption{Computational Study \RN{2}. Comparison of results using the data assimilation algorithm applied to the $2D$ Couette flow with data collected from different regions of the domain: $\Omega_{\text{obst}}$, $\Omega_{\text{bulk}}$, and $\Omega_{\text{wall}}$.}
\label{fig:Couette-Loc}
\end{figure}

Figure~\ref{fig: r1_point2_r2_point9} shows the time evolution of the \( L^2 \) norm of the difference in velocities between the reference solution and the nudged solutions, using both local and global data assimilation. Full synchronization is achieved when data is used from the interior region \( \Omega_{\text{bulk}} \), which lies away from the boundaries, although at a slightly slower rate compared to full-domain assimilation. In contrast, observational data collected near the boundaries (from regions \( \Omega_{\text{obst}} \) and \( \Omega_{\text{wall}} \)) are found to be largely uninformative.

Next, we vary the regions of observed data two more times by adjusting the radius  of the artificial boundaries \( C_1(s) \) and \( C_2(s) \) to
\[
r_1 = 0.3 \quad \text{and} \quad r_2 = 0.8,
\]
and  
\[
r_1 = 0.5 \quad \text{and} \quad r_2 = 0.85.
\]  
\begin{table}[ht!]
\centering
\begin{tabular}{|c c c c|} 
 \hline
 Domain & $\#$ Vertices &    $\#$ Vertices per Area  & Max size of mesh   \\ [0.5ex] 
 \hline
$\Omega$ (DNS) & 10545  & 3401  & 0.032   \\ 
$\Omega_{\text{obst}}$ & 526  & 2140  & 0.043 \\
$\Omega_{\text{bulk}}$ & 1425  &  828  &  0.057 \\
$\Omega_{\text{wall}}$ & 919 & 813  &  0.064  \\ [1ex] 
 \hline
\end{tabular}
\caption{Computational Study \RN{2}. Details on the  Couette spatial meshes for $r_1=0.3, \ r_2=0.8$.}
\label{table:2}
\end{table}

\begin{table}[ht!]
\centering
\begin{tabular}{|c c c c|} 
 \hline
 Domain & $\#$ Vertices &    $\#$ Vertices per Area  & Max size of mesh   \\ [0.5ex] 
 \hline
$\Omega$ (DNS) & 7973  & 2571  & 0.042   \\ 
$\Omega_{\text{obst}}$ & 747  & 996  & 0.071 \\
$\Omega_{\text{bulk}}$ & 733  &  495  &  0.084 \\
$\Omega_{\text{wall}}$ & 747 & 858  &  0.066  \\ [1ex] 
 \hline
\end{tabular}
\caption{Computational Study \RN{2}. Details on the  Couette spatial meshes for $r_1=0.5, \ r_2=0.85$.}
\label{table:3}
\end{table}

In these configurations, the outer subdomains \( \Omega_{\text{obst}} \) and \( \Omega_{\text{wall}} \) become larger, while the interior subdomain \( \Omega_{\text{bulk}} \) becomes smaller; see Figure~\ref{fig: r1_point3_r2_point8} and Figure~\ref{fig: r1_point5_r2_point85}. 

Detailed information about the DNS mesh and the local assimilation meshes is provided in Tables~\ref{table:2} and~\ref{table:3}, respectively. 

The same phenomenon is observed (See Figures \ref{fig: r1_point2_r2_point10} and \ref{fig: r1_point2_r2_point11}): Full synchronization is consistently achieved when using data from the interior region \( \Omega_{\text{bulk}} \), whereas data collected near the boundaries (from regions \( \Omega_{\text{obst}} \) and \( \Omega_{\text{wall}} \)) continue to have minimal to  no impact on the assimilation process. In the final test \ref{fig: r1_point2_r2_point11}, where the interior observation region \( \Omega_{\text{bulk}} \) was considerably smaller, synchronization was still achieved, but only after a substantially longer simulation time. This further illustrates that while interior observations are sufficient, as expected, the size and density of the observational region strongly influence the rate of convergence.

\section{Computational Study III: Forced flow with Non-Homogeneous Dirichlet Boundary Conditions}\label{Section; SIM3}

For our third  experiment,  we consider the flow between offset circles. In short,  comparing to the previous  test in Section \ref{Section; SIM2},  we disrupt the symmetry of the domain and introduce a driven force. The domain consists of a disk with an inner smaller off-center disc removed, defined as
$$
\Omega = \{ (x,y) \in \mathbb{R}^2 : x^2 + y^2 \leq 1^2 \hspace{0.3cm} \text{and} \hspace{0.3cm} (x - 0.2)^2 + y^2 \geq (0.1)^2 \}, 
$$
with no-slip boundary conditions on both circles.  While the Reynolds number  is $\Rey = 600$,  the flow is also driven by a counterclockwise rotational, time-independent body force
$$
\bff (x,y) = (-4y(1 - x^2 - y^2), 4x(1 - x^2 - y^2))^{\tran}.
$$
As the counterclockwise force applies, the flow rotates around $(0, 0)$ and interacts with the immersed circle (obstacle). This induces a von Kármán vortex street which interacts with the near wall streaks common in turbulent flow and a central vortex.

\begin{figure}[ht!]
  \centering
  \begin{subfigure}[b]{0.49\linewidth}
    \includegraphics[width=0.8\linewidth]{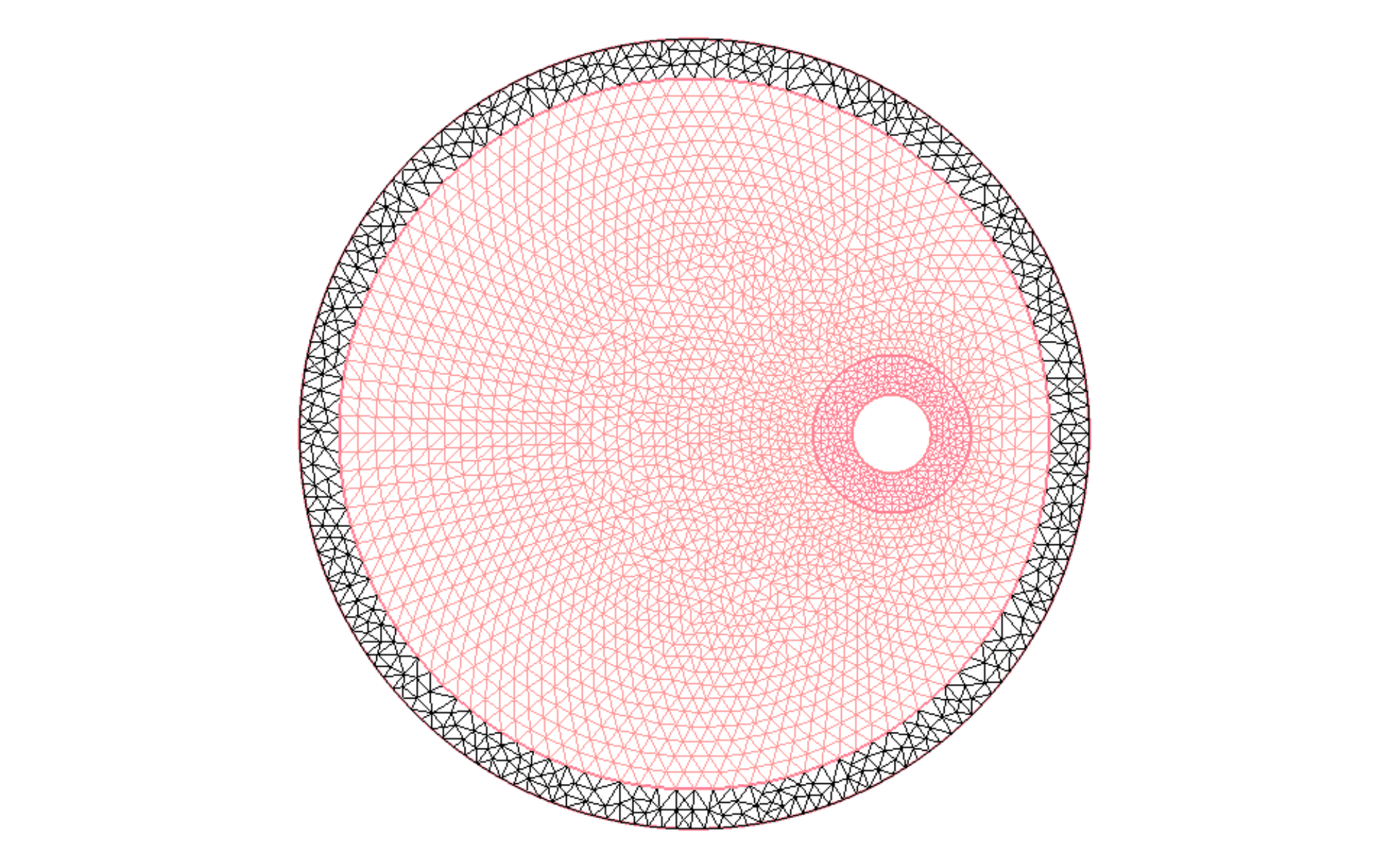}
    \caption{$r_1=0.2,  \ r_2=0.9$.}
  \end{subfigure}
  \begin{subfigure}[b]{0.49\linewidth}
    \includegraphics[width=0.8\linewidth]{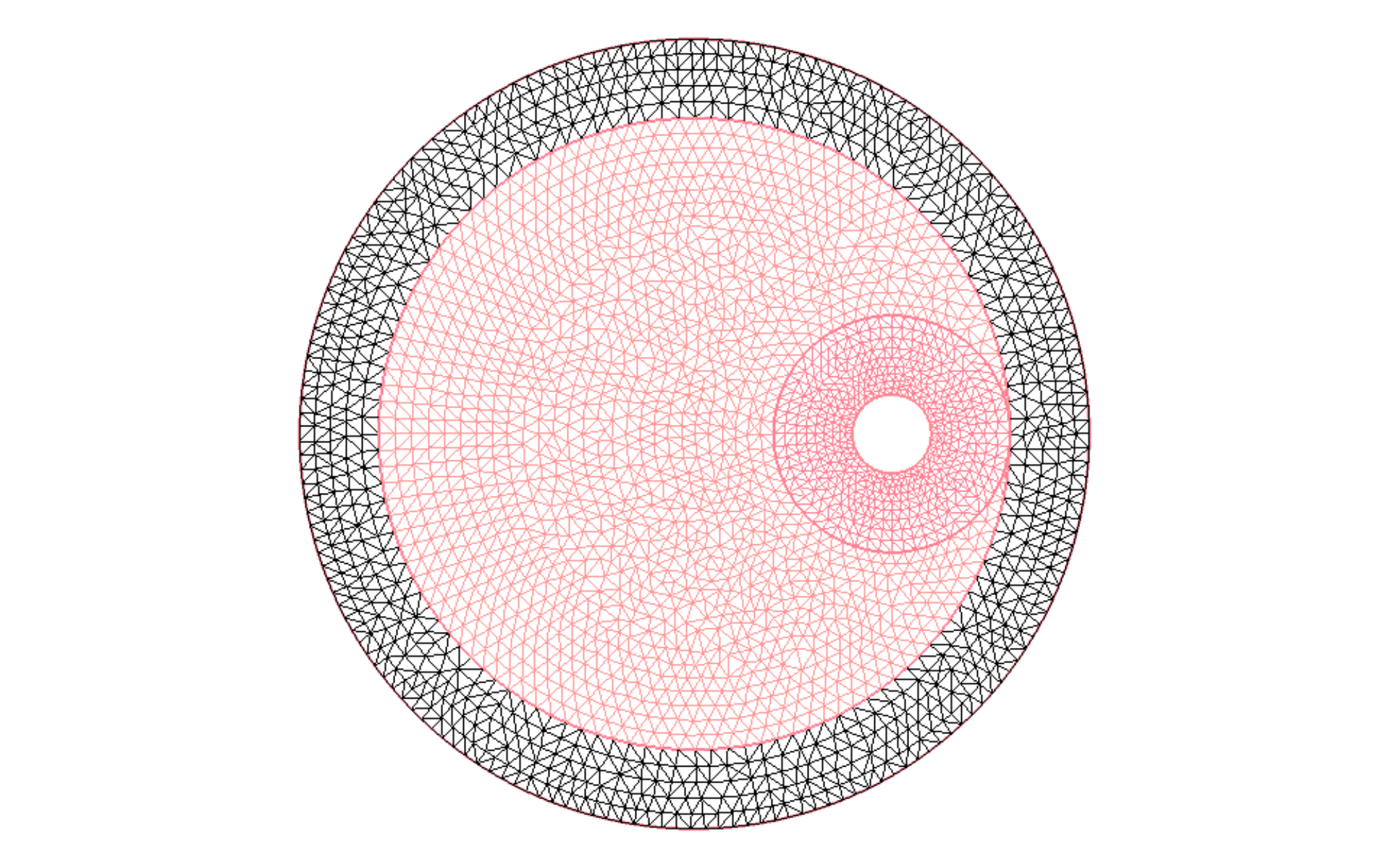}
    \caption{$r_1=0.3,  \ r_2=0.8$.}
  \end{subfigure}
 
  \caption{Computational Study \RN{3}. Local coarse mesh construction within the domain \(\Omega\) showing three sub-domains: \(\Omega_{\text{obst}}\) (closest to the inner disc - dark red), \(\Omega_{\text{bulk}}\) (between the inner and outer regions - light red), and \(\Omega_{\text{wall}}\) (closest to the outer boundary - black).}
  \label{fig:LocMeshBF}
\end{figure}

Using the same approach as before, each sub-domain $\Omega_{\text{obst}}$, $\Omega_{\text{bulk}}$, and $\Omega_{\text{wall}}$ has been adjusted by shifting $0.2$ units towards the x-axis, reflecting the relevant modifications, see Figure~\ref{fig:LocMeshBF}.

The true initial condition $\bfu_0$ is generated by solving the steady Stokes problem with  the same body forces $\bff (x,y)$ and non-homogeneous boundary conditions. Given the absence of an exact solution for this problem, we adopt the DNS solution $\bfu^h$ as our reference solution. 

\subsection*{Impact of Observation Placement: Boundary vs. Interior}
\begin{figure}[ht!]
  \centering
  \begin{subfigure}[b]{0.49\linewidth}
    \includegraphics[width=0.8\linewidth]{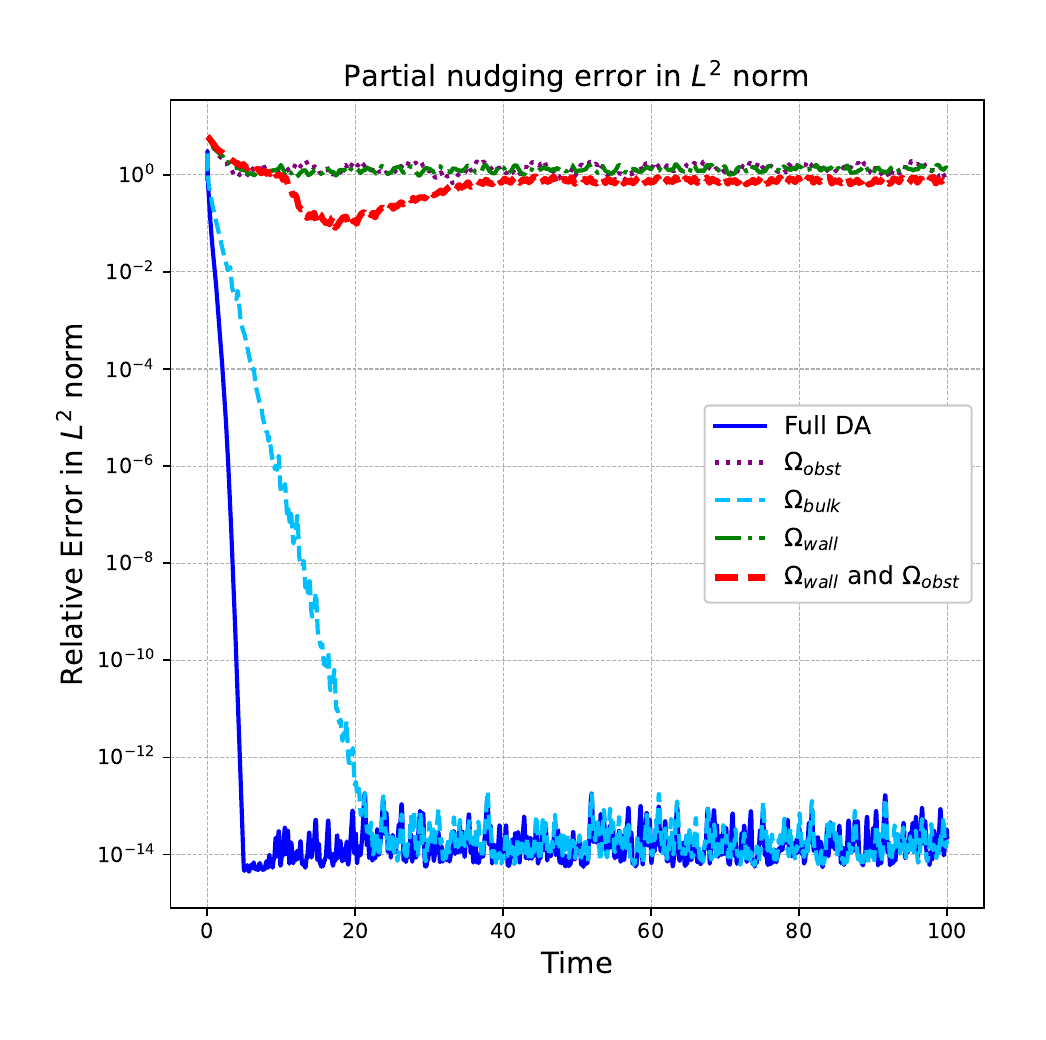}
    \caption{Flow between offset circles $r_1=0.2, \ r_2=0.9$.}
    \label{error_r1_point2_r2_point9}
  \end{subfigure}
  \begin{subfigure}[b]{0.49\linewidth}
    \includegraphics[width=0.8\linewidth]{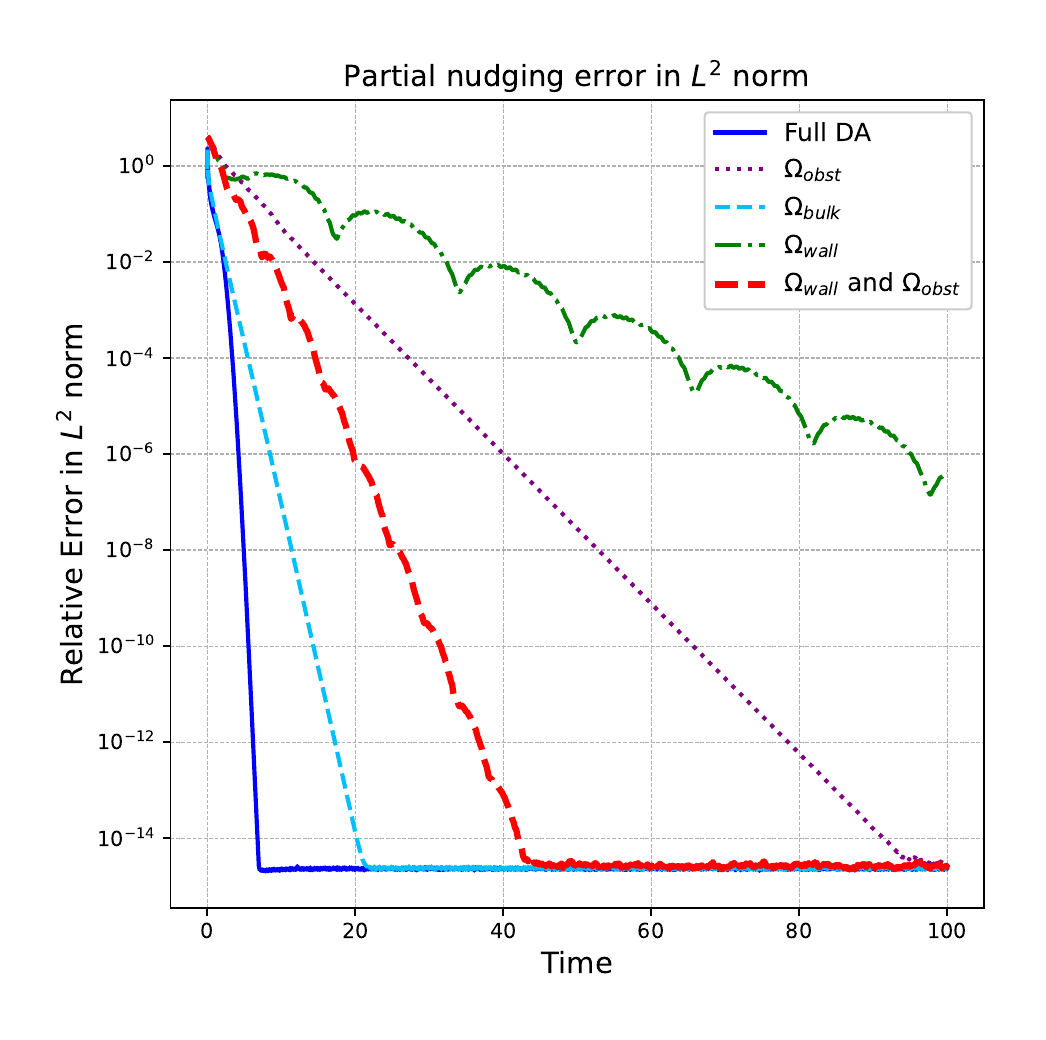}
    \caption{Flow between offset circles $r_1=0.3, \ r_2=0.8$.}
    \label{error_r1_point3_r2_point8}
  \end{subfigure}

\caption{Computational Study \RN{3}. Comparison of results using the DA algorithm applied to the forced flow between offset circles  with data collected from different regions of the domain: $\Omega_{\text{obst}}$, $\Omega_{\text{bulk}}$, and $\Omega_{\text{wall}}$.}

 \label{fig:BF}
\end{figure}

The results from this experiment (see Figures~\ref{error_r1_point2_r2_point9} and \ref{error_r1_point3_r2_point8}) reaffirm the previous findings: assimilation using local data from the interior region \( \Omega_{\text{bulk}} \), away from the boundaries, leads to  global synchronization up to the machine precision, while observations taken near the boundary regions \( \Omega_{\text{obst}} \) and \( \Omega_{\text{wall}} \) provide little to no useful information for the assimilation process.

\FloatBarrier

\section{Discussion}\label{Section; Discussion}

Given the chaotic nature of turbulent flow, estimating the turbulent state by simply interpolating sparse observations is ineffective. DA addresses this challenge by combining observations (collected through sensors or other measurement techniques) with numerical simulations based on the governing physical model, thereby enabling high-fidelity reconstruction of the entire flow state. Due to practical limitations in sensor placement, efficiently reconstructing the full turbulent flow requires strategic positioning of sensors, targeting flow features that most strongly influence the system; so that each observation provides maximal information about the overall dynamics. This challenge is addressed in the present work using the nudging algorithm.  This work highlights the potential of localized DA to reduce observational costs while maintaining predictive fidelity. 

\AAA

In Theorem~\ref{Theorem;Main}, we rigorously prove that global recovery is guaranteed when the observation region is sufficiently close to the entire domain; more precisely, every point in the domain must lie within a distance  
\[
\delta \lesssim H \sim \mathcal{O}(\nu^{1/2}).
\]
from the observed region, and in practice we observe robust synchronization even when this condition is exceeded.

To compare this sufficient condition with the computational regime, we summarize 
in Table~\ref{tab:distance-comparison} the parameters associated with the
interior bulk observation region, $\Omega_{\mathrm{bulk}}$, in Computational
Studies I and II. We also report the dimensionless ratios
\[
R_{\nu}:=\frac{\delta}{\sqrt{\nu}},
\qquad
R_{\delta}:=\frac{\delta}{H},
\qquad\text{and}\qquad
R_H:=\frac{4\mu H^2}{\nu}.
\]

\begin{table}[htb]
\centering

\resizebox{\textwidth}{!}{%
\begin{tabular}{cccccccccc}
\hline
Study & Configuration & $\nu$ & $\mu$ & $H$ & $\delta$
& $R_{\nu}$ & $R_{\delta}$ & $R_H$ & Synchronization \\
\hline
I  & $\ell=1$       & $1$     & $5$  & $1.70$  & $1$
   & $1.00$ & $0.59$ & $57.80$  & Yes \\

I  & $\ell=2$       & $1$     & $5$  & $1.51$  & $2$
   & $2.00$ & $1.33$ & $45.60$  & Yes \\

I  & $\ell=3$       & $1$     & $5$  & $1.38$  & $3$
   & $3.00$ & $2.18$ & $38.09$  & Yes \\

II & $(r_1,r_2)=(0.2,0.9)$
   & $1/600$ & $10$ & $0.064$ & $0.10$
   & $2.45$ & $1.56$ & $98.30$  & Yes \\

II & $(r_1,r_2)=(0.3,0.8)$
   & $1/600$ & $10$ & $0.057$ & $0.20$
   & $4.90$ & $3.51$ & $77.98$  & Yes \\

II & $(r_1,r_2)=(0.5,0.85)$
   & $1/600$ & $10$ & $0.084$ & $0.40$
   & $9.80$ & $4.76$ & $169.34$ & Yes \\
\hline
\end{tabular}%
}
\caption{Comparison of the computational parameters for observations restricted
to $\Omega_{\mathrm{bulk}}$ with the sufficient scaling conditions in
Theorem~\ref{Theorem;Main}. Here,
$R_{\nu}=\delta/\sqrt{\nu}$,
$R_{\delta}=\delta/H$, and
$R_H=4\mu H^2/\nu$.}
\label{tab:distance-comparison}
\end{table}

As shown in Table~\ref{tab:distance-comparison}, all six interior-observation
configurations achieve global synchronization, even though
$\delta/\sqrt{\nu}>1$ and $\delta/H>1$. In particular,
$\delta/\sqrt{\nu}$ ranges from approximately $2$ to $6$ in Computational
Study I and from approximately $2.45$ to $9.80$ in Computational Study II.
Thus, the simulations operate beyond the nominal geometric scaling suggested
by the sufficient analytical condition, while still recovering the global
flow field. This provides quantitative evidence that the theoretical bound is
conservative and can be substantially relaxed in practice. We emphasize that the constants $C_1$ and $C_p$ appearing in
Theorem \ref{Theorem;Main} are not explicitly known. Consequently, the ratios in
Table~\ref{tab:distance-comparison} should not be interpreted as establishing
a strict violation of the theorem's hypotheses. Rather, they provide a
transparent comparison with the nominal order-one scaling underlying the
theoretical estimates.
\PPP

Previous studies employing a range of methodologies~\cites{suzuki2017estimation, bewley2004footprints, encinar2019logarithmic, wang2021state-estimation} have shown that wall data can be used to reconstruct near-wall turbulence and, beyond that, only the outer large-scale structures of the flow. Our results resonate with the findings of Zaki and Wang~\cite{zaki2021flow-reconstruction}, who observed that while wall measurements can reveal near-wall turbulent structures, only the large-scale motions in the outer flow are recoverable from such data. Beyond this, our simulations  confirm that observations collected near the boundary are insufficient to reconstruct the full flow field with high fidelity. However, in contrast to the limitations of wall-based sensing, we demonstrate, both theoretically and through numerical experiments, that it is possible to reconstruct the entire turbulent state, including small-scale structures, by utilizing interior observations located deeper within the bulk of the domain. This suggests that targeted sensing away from boundaries can offer greater predictive power than traditional wall-based approaches, highlighting the potential of strategically placed interior sensors for efficient flow reconstruction in practical settings.

In agreement with Biswas, Bradshaw, and Jolly~\cite{Biswas2021}, we also establish that full flow reconstruction is achievable using observations restricted to a subregion of the domain. While their theoretical results rely on spectral methods on a periodic $2D$ torus, our work extends these insights to a more physically realistic setting: the two-dimensional incompressible NSE on bounded domains with Dirichlet boundary conditions, with simulations performed using finite element methods over complex geometries. Notably, we observe that data near the boundaries are often uninformative.

One notable aspect of our findings is the limited utility of near-boundary observations in achieving global synchronization of the flow. While this may seem counterintuitive at first glance, it is consistent with insights from classical boundary layer theory, which suggest that near-wall dynamics are largely determined by viscosity and boundary conditions.  According to Prandtl's theory, in the high Reynolds number regime, the effects of viscosity are confined to a thin region adjacent to the boundary (known as the boundary layer) whose thickness scales like \( \delta_{\text{BL}} \sim \nu^{1/2} \) \cites{Prandtl1904_BL, Schlichting_Gersten2017_BLTheory}. Within this layer, the velocity field undergoes rapid transition to satisfy the no-slip Dirichlet condition, but its structure is highly constrained: it is slaved to the outer inviscid flow and strongly governed by the viscosity and imposed boundary conditions. Under standard assumptions, such as smooth geometry, regular data, and absence of flow separation, the behavior of the boundary layer is predictable and exhibits a universal profile that depends primarily on the external flow rather than on additional localized perturbations. Consequently, in such regimes, the near-wall dynamics contains limited information; once the viscosity and boundary data are fixed, the wall behavior is essentially determined.

The threshold $\delta \sim \nu^{1/2}$,  in Theorem \ref{Theorem;Main},  mirrors the classical boundary layer thickness, suggesting that the domain outside the layer, where the flow is more turbulent and less constrained, is where uncertainty is concentrated and observations are most valuable. Hence, the success of interior-only nudging can be interpreted as a manifestation of the physical structure of viscous flow: the boundary layer requires no direct correction, because its dynamics are determined indirectly through its coupling with the interior. In this way, our findings not only extend the applicability of continuous data assimilation methods but also reveal a deep and previously unexplored connection between observability in dissipative systems and the classical structure of near-wall flow in fluid dynamics.

An important direction for future research concerns the sharpness of the geometric threshold \( \delta \lesssim \nu^{1/2} \) established in our analysis. While our simulations suggest that synchronization can occur even when this bound is exceeded, it remains an open question whether the \( \nu^{1/2} \) scaling is truly optimal or merely an artifact of the analytical methods employed. Constructing counterexamples in which synchronization provably fails for \( \delta \gg \nu^{1/2} \) would help clarify this issue and further solidify the connection to Prandtl boundary layer theory, where \( \nu^{1/2} \) governs the wall-layer thickness. Another promising avenue is the development of adaptive strategies for selecting the observational region \( \Omega_0 \). Instead of fixing this region a priori, it may be more efficient to place sensors based on dynamically evolving flow features, such as enstrophy, vorticity magnitude, or coherent structures like vortices (see, for example,  the work of Franz, Larios, and Victor \cite{Franz_Larios_Victor2022}). Finally, while our theoretical and computational results are restricted to two-dimensional flows, extending the local DA framework to three-dimensional turbulence presents both mathematical and numerical challenges. In particular, it is unclear whether interior-only nudging remains sufficient for achieving global synchronization in the presence of richer small-scale interactions and more complex boundary layer dynamics.

\AAA
The baseline numerical experiments use noise-free observations generated by a twin DNS model, allowing synchronization to reach machine precision. To assess the effect of observational noise, we also performed simulations with noisy interior data. As shown in Subsection 5.2, the assimilation error converges to a noise-dependent plateau whose level increases with the noise amplitude. We further examined whether supplementing noisy interior observations with wall-adjacent and obstacle-adjacent data improves robustness. For the configuration considered, adding these boundary observations does not significantly reduce the assimilation error, suggesting that they remain largely redundant even in the presence of noisy interior data. Nevertheless, their effectiveness may depend on the noise level, spatial correlations, observation density, and weighting assigned to different observation regions. Moreover, the present experiments do not address structural model error. A broader investigation of correlated noise, alternative observation weighting, and model–observation mismatch is left for future work (see \cites{cibik2025data, DA_Mismatch2026, CIBIKFang2026}).
\PPP

\section*{Acknowledgment}
The authors thank W. J. Layton for helpful  discussions.  RF was partially supported by NSF grant DMS-2410893. AP was partially supported by NSF grant DMS-2532987 and by a 2026–2027 Research, Scholarship and Creative Activity (RSCA) Award from California State University, Northridge.

\bibliographystyle{plain}
\bibliography{mybib}

\end{document}